\documentclass[preprint]{elsarticle}
\makeatletter
\def\ps@pprintTitle{%
 \let\@oddhead\@empty
 \let\@evenhead\@empty
 \def\@oddfoot{\centerline{\thepage}}%
 \let\@evenfoot\@oddfoot}
\makeatother

\usepackage{lineno,hyperref}
\modulolinenumbers[5]

\usepackage[margin=2.5cm]{geometry}

\usepackage{algorithm}
\usepackage[noend]{algpseudocode}
\usepackage{amssymb,amsmath,amsthm}
\usepackage{bbm}
\usepackage[title]{appendix}
\usepackage{bm}
\usepackage{caption}
\usepackage{lipsum}
\usepackage{siunitx}
\usepackage{subfig}
\usepackage{tcolorbox}
\usepackage{xcolor}

\usepackage{algorithm,amsmath,lipsum,xcolor,caption}

\DeclareDocumentCommand \Aij{o}{A_{ij}}
\DeclareDocumentCommand \Ab{o}{\bm{A_{#1}}}

\DeclareDocumentCommand \alphaa{o}{\alpha_{#1}}
\DeclareDocumentCommand \alphab{o}{\bm{\alpha_{#1}}}

\DeclareDocumentCommand \alphabmijx{o}{\bm{\breve{\alpha}}}

\DeclareDocumentCommand \alphah{o}{\hat{\alpha}_{#1}}

\DeclareDocumentCommand \alphabhmijx{o}{\bm{\hat{\alpha}}_{-ij}}

\DeclareDocumentCommand \alphabh{o}{\bm{\hat{\alpha}_{#1}}}
\DeclareDocumentCommand \alphabhu{o}{\bm{\hat{\alpha}^{#1}}}

\DeclareDocumentCommand \Bb{o}{\bm{B_{#1}}}
\DeclareDocumentCommand \Bbt{o}{\bm{\tilde{B}_{#1}}}

\DeclareDocumentCommand \dbh{o}{\bm{\hat{d}_{#1}}}

\DeclareDocumentCommand \eb{o}{\bm{e}_{#1}}

\DeclareDocumentCommand \eps{}{\varepsilon}
\DeclareDocumentCommand \epsb{o}{\bm{\varepsilon}_{#1}}
\DeclareDocumentCommand \epsbi{o}{\bm{\varepsilon}_{#1}}

\DeclareDocumentCommand \epsbI{o}{\bm{\varepsilon_{I}}}

\DeclareDocumentCommand \epsbh{o}{\bm{\hat{\varepsilon}}}
\DeclareDocumentCommand \epsh{o}{\hat{\varepsilon}_{#1}}

\DeclareDocumentCommand \eps{o}{\varepsilon_{#1}}

\DeclareDocumentCommand \feps{}{f_\varepsilon}

\DeclareDocumentCommand \fepsi{o}{f_{\varepsilon_{#1}}}

\DeclareDocumentCommand \fepsik{o o}{{f_{\varepsilon_{#1}}^{(#2)}}}

\DeclareDocumentCommand \fX{}{f_X}

\DeclareDocumentCommand \fY{}{f_Y}
\DeclareDocumentCommand \fhY{}{\hat{f}_Y}

\newcommand{\floor}[1]{\lfloor #1 \rfloor}

\DeclareDocumentCommand \gammab{}{\bm{\gamma}}

\DeclareDocumentCommand \Gth{o}{\tilde{G}_{h,#1}}

\DeclareDocumentCommand \hopt{}{h_\text{opt}}
\DeclareDocumentCommand \hoptRP{}{h_\text{opt}^{RP}}

\DeclareDocumentCommand \Ib{}{\bm{I}}

\DeclareDocumentCommand \mub{o}{\bm{\mu_{#1}}}

\DeclareDocumentCommand \Adj{}{\text{Adj}}

\DeclareDocumentCommand \Bias{}{\text{Bias}}
\DeclareDocumentCommand \Corr{}{\text{Corr}}
\DeclareDocumentCommand \Cov{}{\text{Cov}}
\DeclareDocumentCommand \Det{}{\text{Det}}

\DeclareDocumentCommand \Gamma{}{\text{Gamma}}
\DeclareDocumentCommand \ISE{}{\text{ISE}}
\DeclareDocumentCommand \MISE{}{\text{MISE}}
\DeclareDocumentCommand \MSE{}{\text{MSE}}

\DeclareDocumentCommand \Var{}{\text{Var}}

\DeclareDocumentCommand \mui{o}{\mu_{#1}}
\DeclareDocumentCommand \muN{o}{\mu_{N,#1}}

\DeclareDocumentCommand \oe{o}{\emph{o}}

\DeclareDocumentCommand \Phii{o}{\Phi_{#1}}

\DeclareDocumentCommand \PsiA{o}{\Psi_{#1}^{(A)}}
\DeclareDocumentCommand \PsiAt{o}{\tilde {\Psi}_{#1}^{(A)}}

\DeclareDocumentCommand \Qb{o}{\bm{Q}_{#1}}

\DeclareDocumentCommand \ri{o}{r_{#1}}

\DeclareDocumentCommand \sbi{o}{\bm{s_{#1}}}
\DeclareDocumentCommand \si{o}{s_{#1}}
\DeclareDocumentCommand \st{o}{\tilde{s}_{#1}}

\DeclareDocumentCommand \tbi{o}{\bm{t_{#1}}}
\DeclareDocumentCommand \sigmaeps{}{\sigma}

\DeclareDocumentCommand \vb{o}{\bm{v_{#1}}}
\DeclareDocumentCommand \Vb{o}{\bm{V_{#1}}}

\DeclareDocumentCommand \Wb{o}{\bm{W_{#1}}}

\DeclareDocumentCommand \x{o}{x_{#1}}
\DeclareDocumentCommand \xb{o}{\bm{x_{#1}}}

\DeclareDocumentCommand \X{o}{X_{#1}}
\DeclareDocumentCommand \Xb{o}{\bm{X_{#1}}}
\DeclareDocumentCommand \XbI{o}{\bm{X_{I}}}

\DeclareDocumentCommand \Xbmijx{o}{\bm{X}_{-ij}}

\DeclareDocumentCommand \Xov{o}{\overline{X}_{N,#1}}

\DeclareDocumentCommand \y{o}{y_{#1}}
\DeclareDocumentCommand \yt{}{\tilde{y}}
\DeclareDocumentCommand \yti{o}{\tilde{y}_{#1}}
\DeclareDocumentCommand \ybmijx{o}{\bm{y}_{-ij}}

\DeclareDocumentCommand \Y{o}{Y_{#1}}
\DeclareDocumentCommand \Yb{}{\bm{Y}}

\DeclareDocumentCommand \zi{o}{z_{#1}}

\DeclareDocumentCommand \Z{o}{Z_{#1}}
\DeclareDocumentCommand \Zb{o}{\bm{Z_{#1}}}

\DeclareDocumentCommand \zerob{}{\bm{0}}

\DeclareDocumentCommand \zetab{o}{\bm{\zeta_{#1}}}

\DeclareDocumentCommand \zerob{o}{\bm{0}}

\DeclareDocumentCommand \indOLS{o}{\mathbbm{1}_{\{N^{1/2}|(\alphabh[] - \alphabh[I])_{j}| > #1\}}}

\newtheorem{theorem}{Theorem}[section]
\newtheorem{corollary}[theorem]{Corollary}
\newtheorem{lemma}[theorem]{Lemma}
\newtheorem{proposition}[theorem]{Proposition}

\newtheoremstyle{claim}
  {\topsep}
  {\topsep}
  {\itshape}
  {}
  {}
  {.}
  {.5em}
  {\thmname{#1}\thmnumber{ #2}\thmnote{ (#3)}}
\theoremstyle{claim}

\newtheorem{asu}{Assumption}

\makeatletter
\def\algbackskip{\hskip-\ALG@thistlm}
\makeatother

\begin{document}

\begin{frontmatter}


\title{A Multiple Regression-Enhanced Convolution Estimator for the Density of a Response Variable in the Presence of Additional Covariate Information}

\author{Brian Fitzpatrick\corref{mycorrespondingauthor}}
\cortext[mycorrespondingauthor]{Corresponding author}
\ead{brian.fitzpatrick@tudublin.ie}

\author{James Loughman}
\ead{james.loughman@tudublin.ie}

\author{Daniel Ian Flitcroft}
\ead{ian.flitcroft@tudublin.ie}

\address{Centre for Eye Research Ireland, School of Physics \& Clinical \& Optometric Sciences, TU Dublin, Ireland}

\begin{abstract}
In this paper we propose a convolution estimator for estimating the density of a response variable that employs an underlying multiple regression framework to enhance the accuracy of density estimates through the incorporation of auxiliary information. Suppose we have a sample consisting of $N$ complete case observations of a response variable and an associated set of covariates, along with an additional sample consisting of $M$ observations of the covariates only. We show that the mean square error of the multiple regression-enhanced convolution estimator converges as $O(N^{-1})$ towards zero, and moreover, for a large fixed $N$, that the mean square error converges as $O(M^{-4/5})$ towards an $O(N^{-1})$ constant. This is the first time that the convergence of a convolution estimator with respect to the amount of additional covariate information has been established.
In contrast to convolution estimators based on the Nadaraya-Watson estimator for a nonlinear regression model, the multiple regression-enhanced convolution estimator proposed in this paper does not suffer from the curse of dimensionality. It is particularly useful for scenarios in which one wants to estimate the density of a response variable that is challenging to measure, while being in possession of a large amount of additional covariate information. In fact, an application of this type from the field of ophthalmology motivated our work in this paper. 
\end{abstract}


\end{frontmatter}

\nolinenumbers

\section{Introduction}
The standard approach to estimating the unknown probability density function of a random variable $Y$ is kernel density estimation, a nonparametric statistical technique which can be traced back to the pioneering works of Rosenblatt \cite{rosenblatt1956} and Parzen \cite{parzen1962estimation} over fifty years ago. Conventional kernel density estimation involves estimating the density $\fY$ of $Y$ using the Rosenblatt–Parzen density estimator
\begin{align} \label{eq:kde}
\fhY(y) = \frac{1}{hN} \sum_{i=1}^N K_h(y - Y_i),
\end{align}
where the set $\{Y_i\}_{i=1}^N$ is a sample of $N$ i.i.d observations of $Y$, $K_h(\cdot) = K(h^{-1}(\cdot))$ with $K$ being some kernel function, and $h>0$ is the bandwidth. Recently, there has been a lot of interest in another type of density estimator known as a convolution estimator. A convolution estimator can be employed when $Y$ is related to a set of covariates through a regression model such as
\begin{align} \label{eq:reg-nonlin}
Y = m(X) + \varepsilon,
\end{align}
where $m$ is a regression function, the covariate vector $X$ and the error $\varepsilon$ are independent, and $\varepsilon$ has mean zero and finite variance. The naming convention arises due to the fact that the probability distribution of a summation of random variables can be expressed in terms of a convolution. Estimating the density $\fY$ of $Y$ with a convolution estimator involves first estimating the underlying regression function $m$.

Escancianoa and Jacho-Ch{\'a}vez \cite{escanciano2012n} used the Nadaraya–Watson estimator to estimate the underlying regression function, and established asymptotic normality of their convolution estimator. M\"uller \cite{muller2012estimating} approached the problem in terms of an arbitrary estimator for the underlying regression function, and showed that the convolution estimator can achieve the optimal parametric convergence rate $\sqrt{N}$. St{\o}ve and Tj{\o}stheim \cite{stove2012convolution}, who also employed the Nadaraya–Watson estimator for the underlying regression function, derived explicit expressions for the asymptotic bias and variance of their convolution estimator, and proved that the mean square error (MSE) converges as $O(N^{-1})$.

Li and Tu \cite{li2016n} estimated the underlying regression function using nonlinear least squares, and investigated  important topics such as endogeneity and robustness to misspecification in the regression function, along with proving the $\sqrt{N}$-consistency and asymptotic normality of their convolution estimator.

It is also worth mentioning that both St{\o}ve and Tj{\o}stheim, and  Li and Tu, considered the case when the error can be heteroskedastic. Some other relatively recent works featuring convolution estimators are \cite{schick2004root,schick2007root,saavedra1999rate,saavedra2000estimation}.

By exploiting special structure of $Y$ in \eqref{eq:reg-nonlin}, convolution estimators can achieve $\sqrt{N}$-consistency, and thereby converge much faster than the conventional kernel density estimator which is only $\sqrt{Nh}$-consistent. That said, the convergence of convolution estimators is tied to the convergence of the estimator for the underlying regression model; Muller \cite{muller2012estimating} showed that $\sqrt{N}$-consistency requires plugging in an efficient regression function estimator, while St{\o}ve and Tj{\o}stheim found that their Nadaraya-Waton-based convolution estimator suffers from the curse of dimensionality \cite{stove2012convolution,li2016n}.

A powerful feature of convolution estimators, which has been largely unexplored in existing works on this topic, is that they provide a convenient mechanism by which additional covariate observations, above and beyond the covariate observations in the complete case dataset (response variable and associated covariates) used to define the regression model \eqref{eq:reg-nonlin}, can be incorporated into the estimation process in a straightforward fashion. Denote by $N$ the number of observations in the complete case sample, and by $M$ the number of observations in an additional sample featuring the covariates only. In all of the aforementioned works, apart from M\"uller \cite{muller2012estimating}, the total number of covariate observations matches the number of observations of the response variable, that is, $M=0$. However, there is no reason why the total number of covariate observations can't be larger than the number of observations of the response variable, that is, $M > 0$. This is very interesting as it raises the possibility of enhancing density estimates of a response variable without needing more response observations; instead, additional observations of the covariates can be used to enhance the density estimates.

M\"uller \cite{muller2012estimating} investigated a scenario involving a dataset in which some of the response observations are missing at random, while all of the covariate observations are present. Another interpretation of this situation is that one is in possession of a complete case sample featuring $N$ observations of a response variable and an associated set of covariates, along with an additional sample featuring $M$ observations of the covariates only. This is the perspective we take in this paper.

While the convergence of convolution estimators with respect to $N$ has been established as discussed above, the convergence with respect to $M$ is an open question. This is an important question; it would be useful to know just how effective the incorporation of additional covariate observations is in terms of enhancing the accuracy of density estimates, since often-times in practical applications it can be difficult if not downright impossible to obtain more observations of a response variable, while at the same time it can be very straightforward to obtain more observations of the covariates. For instance, we may want to estimate the density of a response variable that is difficult to measure due to time and/or cost constraints. Since this variable is challenging to measure, it is quite possible that only a small sample of measurements is available. On the other hand, we may find it easy to take or obtain a large number of measurements of other variables that are correlated with the difficult to measure response. To improve the accuracy of estimates of the density of the difficult to measure response, one can incorporate the abundant auxiliary covariate information using a convolution estimator.

In fact, our work in this paper was inspired by an application of this type from the field of ophthalmology. Measurement of the axial length of the human eye has historically been confined to specialist practice areas of ophthalmology, most notably for cataract and refractive surgery. It has not been measured routinely beyond this due to the high cost of biometric devices which are capable of measuring axial length precisely. Axial length has recently emerged as the most important clinical parameter required for the medical management of myopia, a condition associated with excessive eye growth and consequential ocular tissue damage and disease. New treatments are available to limit eye growth in children at risk of progressive myopia, but the clinicians tasked with prescribing and monitoring the efficacy of such treatments do not typically have access to the expensive specialised biometry devices. Consequently, accurate estimates of the axial length distribution in human populations are required to better understand, treat and monitor this and other ocular diseases. Datasets featuring axial length information are limited and small, whereas datsets orders of magnitude greater in size featuring measurements of ocular parameters such as refractive error, corneal radius, and age are readily available.


The convolution estimator proposed in this paper is based on the ordinary least squares estimator (OLS) estimator for an underlying multiple regression framework. This is in contrast to previous works on convolution estimators which have generally considered nonlinear regression functions and utilized nonlinear estimators such as the Nadaraya-Watson estimator. Nononparametric methods such as the Nadaraya-Watson estimator are afflicted by the curse of dimensionality; their convergence scales badly as the dimensionality of the covariates increases \cite{gyorfi2006distribution}. The authors of \cite{stove2012convolution} note that their Nadaraya-Watson based estimator is not suitable for covariate vectors with more than three dimensions for this reason. On the other hand, the convergence of the OLS estimator for a multiple regression model is independent of the number of covariates, since multiple regression is an additive model. Therefore, it is reasonable to expect that the curse of dimensionality will not be an issue for a convolution estimator that employs the OLS estimator for the underlying regression model. Note that by 'multiple regression', we mean a regression model that is linear in the parameters but potentially non-linear in the covariates, such as polynomial regression.


The key issues to consider when deciding on an underlying regression framework for a convolution estimator are (i) the level of nonlinearity present in the data, and (ii) the dimensionality of the covariates. If the data is highly nonlinear with low-dimensional covariate vectors, a Nadaraya-Watson-based estimator for a nonlinear regression model is a strong choice. On the other hand, if the data can be be well fit by a linear model, possibly after some non-linear transformations, or if one wants to use covariate vectors that span many dimensions, the OLS estimator and multiple regression may be a better choice.

Another aspect of convolution density estimators worth highlighting is that they are considerably more computationally expensive than the Rosenblatt–Parzen density estimator, since an evaluation with a convolution estimator requires two summations over the sample data instead of one, for each point on the evaluation grid. This can lead to high computational costs, so techniques for accelerating the computation of convolution estimator evaluations are desirable.

Our focus in this work is on establishing the theoretical and computational foundations of the multiple regression-enhanced convolution estimator. Our applied work on the estimating the distribution of the axial length of the human eye using this estimator will be reported in a future ophthalmology-focused research article. The main contributions of this work are as follows.
\begin{enumerate}
\item We derive the asymptotically optimal bandwidth for the multiple regression-enhanced convolution estimator, and show that it can be related to the asymptotically optimal bandwidth for the classical Rosenblatt–Parzen density estimator. In particular, the dependence of the optimal bandwidth on both $N$ and $M$ is established.
\item We show that the MSE of the multiple regression-enhanced convolution estimator converges as $O(N^{-1})$ irrespective of the dimensionality of the covariates, which means that it is not afflicted by the curse of dimensionality.
\item We resolve the question on the convergence of convolution estimators with respect to the number of covariate observations in the additional sample, by showing that for a large fixed $N$, the MSE converges as $O(M^{-4/5})$ towards an $O(N^{-1})$ constant. In other words, the accuracy improvement achievable through the incorporation of additional covariate observations eventually saturates at a level that is dependent on the number of complete case samples used in the underlying multiple regression model.
\item We develop a Fast Gauss Transform-based algorithm that substantially reduces the amount of computational time needed to perform convolution density estimator evaluations.
\end{enumerate}

This paper is structured as follows. In Section \ref{sec:MRBCDE}, we define the multiple regression-enhanced convolution density estimator and state some assumptions that are necessary for the mathematical analysis of the estimator, while also introducing some notational conventions.

In Section \ref{sec:theoretical-analysis}, we present our theoretical analysis which involves deriving the asymptotic bias and variance of the convolution estimator. 

In Section \ref{sec:bandwidth-selection}, we derive the asymptotically optimal bandwidth for the convolution estimator, in particular showing how it depends on both $N$ and $M$. Moreover, we derive the rate of convergence of the MSE of the convolution estimator with respect to both $N$ ad $M$.

In Section \ref{sec:numerical-implementation}, we consider numerical implementation of the convolution estimator. We propose a computational algorithm that incorporates the high-performance C++ library FIGTree \cite{morariu2008automatic}. This library combines the (Improved) Fast Gauss Transform \cite{greengard1991fast} and Approximate Nearest Neighbor searching \cite{arya1993approximate} to reduce the computational complexity of Gauss transform evaluations.

In Section \ref{sec:numerical-simulations}, we perform a series of numerical simulations to gain an understanding of the convolution estimator's performance and investigate the potential reduction in MISE through the incorporation of additional covariate observations.

The paper ends with some concluding remarks in Section \ref{sec:concluding-remarks}. Appendix \ref{appendix:asy-expectations} features the asymptotic analysis of expectations that arise during the derivation of the asymptotic bias and variance. Appendix \ref{sec:appendix-E-PhiiN-inv-sqr-ord-mag} contains some technical proofs that are required to establish the order of magnitude of a specific term that arises in the bias and variance.

\section{Multiple regression-enhanced convolution density estimator} \label{sec:MRBCDE}
Without loss of generality, we assume that the multiple regression model that we are interested in, which is linear in the parameters but potentially nonlinear in the covariates, has if necessary been converted to a multiple linear regression model by variable transformations. Thus, let $\{(\Y[i], \Xb[i])\}_{i=1}^N$ be a sample of $N$ i.i.d. complete case observations of a random vector $(Y,X)$, where $Y$ is related to the $J$-dimensional covariate vector $X$ through the following multiple regression model
\begin{align} \label{eq:pop-mult-reg}
Y = X^T \alphab[] + \varepsilon.
\end{align}
Here, $\alphab[] = [\alphaa[0],\alphaa[1],\dots,\alphaa[J]]^T$, with $\alphaa[0] \neq 0$ for $i \in \{0,1,\dots,J\}$, is the vector of regression coefficients, and the first element of the covariate vector $X = [1,X_1,X_2,\dots,X_J]^T$ is defined to be one for convenience. The assumptions on the error $\varepsilon$ will be specified later. We are interested in estimating the probability density function $\fY$ of $Y$.

Let $\{\Xb[i]\}_{i=N+1}^L$, where $L = N+M$, be an additional sample of $M$ i.i.d. observations of the covariate vector $X$ only. While we could estimate $\fY$ directly using kernel density estimation applied to the $N$ observations of $Y$, instead we will leverage both the regression model \eqref{eq:pop-mult-reg} and the full set of $L$ covariate observations to provide more accurate density estimates than those given by the conventional approach.

The multiple regression model associated with the complete case dataset $\{(\Y[i], \Xb[i])\}_{i=1}^N$ is
\begin{align} \label{eq:samnple-mult-reg}
\Yb = \Xb[]\alphab[] + \epsb[],
\end{align}
where
\begin{align*}
\Yb
=
\begin{bmatrix}
\Y[1] \\
\Y[2] \\
\vdots \\
\Y[N]
\end{bmatrix},
\quad 
\Xb[]
=
\begin{bmatrix}
\Xb[1]^T \\
\Xb[2]^T \\
\vdots \\
\Xb[N]^T
\end{bmatrix}
=
\begin{bmatrix}
1 & \X[11] & \dots & \X[1J] \\
1 & \X[21] & \dots & \X[2J] \\
\vdots & \vdots & \ddots & \vdots \\
1 & \X[N1] & \dots & \X[NJ] \\
\end{bmatrix}
\quad
\alphab[]
=
\begin{bmatrix}
\alpha_0 \\
\alpha_1 \\
\vdots \\
\alpha_J
\end{bmatrix},
\quad
\epsb[]
=
\begin{bmatrix}
\eps[1] \\
\eps[2] \\
\vdots \\
\eps[N]
\end{bmatrix}.
\end{align*}
Denote by $\alphabh[] = [\alphah[0],\alphah[1],\dots,\alphah[J]]^T$ the OLS estimator for the coefficient vector $\alphab[]$. The OLS estimator is given by \cite[4.4]{greene2003econometric}
\begin{align} \label{eq:OLS-estimator}
\alphabh[] = \alphab[] + (\Xb[]^T \Xb[])^{-1} \Xb[]^T \epsb[].
\end{align}
Denote by
\begin{align} \label{eq:Phi}
\Phii[N] := N^{-1} \Xb[]^T \Xb[],
\end{align}
and note that since $\Xb[]^T \epsb[] = \sum_{i=1}^N \Xb[i]\eps[i]$, the OLS estimator can be expressed as
\begin{align} \label{eq:OLS-estimator-Phi}
\alphabh[] = \alphab[] + N^{-1} \Phii[N]^{-1} \sum_{i=1}^N \Xb[i]\eps[i].
\end{align}
The residual vector is $\epsbh = \Yb - \Xb[]\alphabh[] = [\epsh[0],\epsh[1],\dots,\epsh[N]]^T$. Since $\Y[]$ is the sum of random variables, its density can be written as a convolution. Denoting by $\feps$ the error density, and by $F$ the covariate distribution, it holds that \cite{muller2012estimating}
\begin{align} \label{eq:fY-conv-form}
\fY(y) = \int \feps(y - \xb[]^T \alphab[]) F(d\xb[]) = E[\feps(y - X^T \alphab[])].
\end{align}
More generally, it holds that
\begin{align} \label{eq:fYk-conv-form}
\fY^{(k)}(y) = E[\feps^{(k)}(y - X^T \alphab[])],
\end{align}
where $\fY^{(k)}$ is the $k$-th derivative of $\fY$. The OLS estimator and the full set of $L$ covariate observations can be used to estimate the right hand side of \eqref{eq:fY-conv-form}:
\begin{align} \label{eq:fY-approx}
\fY(y) \approx \frac{1}{L} \sum_{i=1}^L \feps(y - \Xb[i]^T \alphabh[]).
\end{align}
Next, the residuals and conventional kernel density estimation can be used to estimate $\feps(y)$:
\begin{align} \label{eq:feps-approx}
\feps(y) \approx \frac{1}{N} \sum_{i=1}^N K_h(y-\epsh[i]),
\end{align}
Using \eqref{eq:fY-approx} and \eqref{eq:feps-approx}, we define the multiple regression-enhanced convolution estimator $\fhY$ by
\begin{align} \label{eq:fhY}
\fhY(y) = \frac{1}{hNL} \sum_{i=1}^L \sum_{j=1}^N K_h(y - \Xb[i]^T \alphabh[]-\epsh[j]).
\end{align}
This is the form of the estimator we use for computation. For the mathematical analysis, it is convenient to work with a slightly different expression for the estimator. Noting that $\epsh[j] = \y[j] - \Xb[j]^T \alphabh[]$, it is straightforward to show that $\fhY$ can be written as
\begin{align} \label{eq:fhY-convenient}
\fhY(y)
&= \frac{1}{hNL} \sum_{i=1}^L \sum_{j=1}^N K_h(y - \Xb[i]^T \alphab[] - \eps[j] + (\Xb[j] - \Xb[i])^T(\alphabh[]  - \alphab[])).
\end{align}

%
%
%
%
%
%
%
%
\subsection{Notation} \label{subsec:Notation}
We introduce a function $\yt$ for notational convenience:
\begin{align} \label{eq:yt}
\yt(\gammab,\Xb[i],\Xb[j]) &:= y - \Xb[i]^T \alphab[] + (\Xb[j] - \Xb[i])^T(\gammab - \alphab[])
\end{align}
Note that $\yt$ reduces to a particularly simple form in certain cases, that is,
\begin{align} \label{eq:yt-simple}
\yt(\gammab,\Xb[i],\Xb[j]) = y - \Xb[i]^T \alphab[], \quad \quad \text{for} \ \gammab = \alphab[], \ \text{or} \ \Xb[i] = \Xb[j].
\end{align}
Using \eqref{eq:yt}, the convolution estimator \eqref{eq:fhY-convenient} can be written as
\begin{align} \label{eq:fhY-rewritten}
\fhY(y)
&= \frac{1}{hNL} \sum_{i=1}^L \sum_{j=1}^N K_h(\yt(\alphabh[],\Xb[i],\Xb[j]) - \eps[j]).
\end{align}

While we have explicitly defined the elements of the random vectors encountered above, for convenience we denote by $(\vb[])_i$ the $i$-th element of a random vector $\vb[]$, since this makes it easier to work with more complicated random vectors. Similarly, we denote by $(\Ab[])_{ij}$ the $(i,j)$-th element of a random matrix $\Ab[]$. By an abuse of notation, since $\vb[]$ is random not deterministic, we write $\partial f(\vb[])/\partial(\vb[])_{i}$ for the partial derivative of a function $f$ with respect to the $i$-th element of its vector-valued argument.

We make the following definitions for convenience.
\begin{align}
K_{ij}(\gammab) &:= K_h(\yt(\gammab,\Xb[i],\Xb[j])-\eps[j]), \label{eq:K-ij} \\
\Cov_{ijkl}(\gammab) & := \Cov[K_h(\yt(\gammab,\Xb[i],\Xb[j]) - \eps[j]),K_h(\yt(\gammab,\Xb[k],\Xb[l]) - \eps[l])]. \label{eq:Cov-ijkl}
\end{align}
Also, we denote by
\begin{align} \label{eq:int-K-defs}
\mu_K := \int r^2 K(r) dr, \quad \quad \sigma_K := \int K^2(r) dr, \quad \quad \sigma_{K,2} := \int r^2 K^2(r) dr.
\end{align}

\subsection{Assumptions} \label{subsec:Assumptions}
\begin{asu} \label{asu:K}
$K(y) = (2\pi)^{-1/2} e^{-\frac{1}{2}y^2}$.
\end{asu}
\begin{asu} \label{asu:second-moments}
$E[Y^2] < \infty$, and $E[X_i^2] < \infty$, for $i=1,\dots,J$.
\end{asu}
\begin{asu} \label{asu:iid}
The $N$ observations in dataset $\{(\Y[i],\Xb[i])\}_{i=1}^N$, and the $M$ observations in the dataset $\{\Xb[i]\}_{i=N+1}^L$ are independent and identically distributed.
\end{asu}
\begin{asu} \label{asu:error}
$E[\epsb[]|\Xb[]] = 0$, $E[\epsb[i]^2|\Xb[]] = \sigmaeps^2$, and $E[\epsb[i]\epsb[j]|\Xb[]] = 0$ for $i \neq j$.
\end{asu}
\begin{asu} \label{asu:diff}
The error density $\feps$ is four times differentiable.
\end{asu}
\begin{asu} \label{asu:n-sqr-h}
The bandwidth $h = h(N)$ behaves as $\lim_{N \to \infty} h(N) = 0$, and $\lim_{N \to \infty} h(N)N^2 = \infty$.
\end{asu}

Assumption \ref{asu:K} means that we are restricting to the Gaussian kernel function. We have restricted the kernel to the Gaussian function because our computational implementation of the multiple regression-enhanced convolution estimator is based on the Fast Gauss Transform. Due to Assumption \ref{asu:K}, all derivatives of $K$ are bounded. Moreover,
\begin{align} \label{eq:int-K}
\int K(r) dr = \int r^2 K(r) dr = 1, \quad \quad \int r K(r) dr = 0, \quad \quad \int r K^2(r) dr = 0.
\end{align}
Assumptions \ref{asu:second-moments}, \ref{asu:iid}, and \ref{asu:error} are standard conditions for multiple regression. Assumption \ref{asu:diff} is a standard regularity condition that ensures well-defined Taylor expansions. Assumption \ref{asu:n-sqr-h} is analogous to the usual assumption in conventional kernel density estimation that ensures the variance converges to zero as $N \to \infty$.
\section{Theoretical Analysis} \label{sec:theoretical-analysis}
%
%
%
%
We begin by deriving the asymptotic bias and variance of the multiple regression-enhanced convolution estimator.

\begin{theorem}[] \label{thm:bias-asy}
Under assumptions (A), (B), (C), (D), (E), and (F), the asymptotic bias of the multiple regression-enhanced convolution estimator is
\begin{align} \label{eq:bias}
\Bias[\fhY(y)]
& = O(h^2) + O\bigg(\frac{1}{N}\bigg).
\end{align}
\begin{proof}
First,
\begin{align*}
E[\fhY(y)]
&= \frac{1}{hNL} \sum_{i=1}^L \sum_{j=1}^N E[K_h(\yt(\alphabh[],\Xb[i],\Xb[j]) - \eps[j])] \\
&= \frac{1}{hNL} \bigg(\sum_{i=1}^N \sum_{j=1}^N + \sum_{i=1}^L \sum_{\substack{j=1 \\ j\neq i}}^N\bigg) E[K_h(\yt(\alphabh[],\Xb[i],\Xb[j]) - \eps[j])] \\
&= \frac{N}{hNL} E[K_h(\yt(\alphabh[],\Xb[1],\Xb[1]) - \eps[1])]
+ \frac{N(N-1)+ MN}{hNL} E[K_h(\yt(\alphabh[],\Xb[1],\Xb[2]) - \eps[2])] \\
&= \frac{1}{hL} E[K_h(y - \Xb[1]^T \alphab[] - \eps[1])]
+ \frac{L-1}{hL} E[K_h(\yt(\alphabh[],\Xb[1],\Xb[2]) - \eps[2])],
\end{align*}
where we used \eqref{eq:yt-simple} for the first expression on the last line.
Then, by Lemma \ref{lem:E-K-11-leading-order},
\begin{align*}
\frac{1}{hL} E[K_h(y - \Xb[1]^T \alphab[] - \eps[1])]
\sim \frac{1}{hL} (h \fY(y) + h^3 \frac{\mu_K}{2} \fY''(y))
= \frac{1}{L} (\fY(y) + h^2 \frac{\mu_K}{2} \fY''(y)).
\end{align*}
Next, by Lemma \ref{lem:E-K-12}, as $N \to \infty$,
\begin{align*}
\frac{L-1}{hL} E[K_h(\yt(\alphabh[],\Xb[1],\Xb[2]) - \eps[2])]
& \sim \frac{L-1}{L}\bigg(\fY(y) + h^2 \frac{\mu_K}{2} \fY''(y) + N^{-1} \frac{\sigmaeps^2}{2} \sum_{p_1,p_2=0}^J E[(\Phii[N]^{-1} \Xb[1])_{p_1}(\Phii[N]^{-1} \Xb[1])_{p_2}] \\
& \times E[f_{\eps[]}''(y - \Xb[1]^T \alphab[]) (\Xb[2] - \Xb[1])_{p_1}(\Xb[2] - \Xb[1])_{p_2}]\bigg).
\end{align*}
Combining these results, we have that
\begin{align*}
E[\fhY(y)]
& \sim \fY(y) + h^2 \frac{\mu_K}{2} \fY''(y)
+ N^{-1} \frac{\sigmaeps^2}{2} \sum_{p_1,p_2=0}^J E[(\Phii[N]^{-1} \Xb[1])_{p_1}(\Phii[N]^{-1} \Xb[1])_{p_2}] \\
& \times E[f_{\eps[]}''(y - \Xb[1]^T \alphab[]) (\Xb[2] - \Xb[1])_{p_1}(\Xb[2] - \Xb[1])_{p_2}].
\end{align*}
Now, $E[(\Phii[N]^{-1} \Xb[1])_{p_1}(\Phii[N]^{-1} \Xb[1])_{p_2}] = O(1)$ as $N \to \infty$ by Lemma \ref{lem:E-prod-PhiiN-Xi}. Therefore, the result follows since $\Bias[\fhY(y)] = E[\fhY(y)] - \fY(y)$.
\end{proof}
\end{theorem}


Before we derive the variance, we need a couple of lemmas. These lemmas, which hold under assumptions (A), (B), (C), (D), (E), and (F) given in Section \ref{subsec:Assumptions}, provide leading-order expressions for terms that arise when we perform a decomposition of the variance in Theorem \ref{thm:var-asy}. It transpires that only the terms $\Cov_{ijkl}(\alphabh[])$ for $(i,j,k,l) \in \{(1,2,1,2),(1,2,1,3),(1,2,3,2),(1,2,3,4)\}$ are important asymptotically. Expressions for these terms are derived in  Appendix \ref{appendix:asy-expectations}. The remaining $\Cov_{ijkl}(\alphabh[])$ terms can be handled in a similar fashion so we omit the repetitive derivations. 

%

%
%
%
%
\begin{lemma} \label{lem:Cov-1212}
It holds that
\begin{align*}
\Cov_{1212}(\alphabh[])
\sim O(h).
\end{align*}
\begin{proof}
Since $\Cov_{1212}(\alphabh[]) = E[K_{12}^2(\alphabh[])] - E[K_{12}(\alphabh[])]^2$, by Lemma \ref{lem:E-K-12-sqrt-internal} and Corollary \ref{cor:E-K-12-sqr}, it holds that
\begin{align*}
\Cov_{1212}(\alphabh[])
\sim h \sigma_{K} \fY(y).
\end{align*}
\end{proof}
\end{lemma}

%
%
%
%
\begin{lemma} \label{lem:Cov-1213}
It holds that
\begin{align*}
\Cov_{1213}(\alphabh[])
\sim O(h^2).
\end{align*}
\begin{proof}
Since $\Cov_{1213}(\alphabh[]) = E[K_{12}(\alphabh[])K_{13}(\alphabh[])] - E[K_{12}(\alphabh[])]^2$, by Lemma \ref{lem:E-K-1213} and Corollary \ref{cor:E-K-12-sqr},
\begin{align*}
\Cov_{1213}(\alphabh[])
\sim h^2 (E[f_{\eps[]}^2(y - \Xb[1]\alphab[])] - \fY^2(y)).
\end{align*}
\end{proof}
\end{lemma}

%
%
%
%
\begin{lemma} \label{lem:Cov-1232}
It holds that
\begin{align*}
\Cov_{1232}(\alphabh[])
\sim O(h^2).
\end{align*}
\begin{proof}
Since $\Cov_{1232}(\alphabh[]) = E[K_{12}(\alphabh[])K_{32}(\alphabh[])] - E[K_{12}(\alphabh[])]^2$, by Lemma \ref{lem:E-K-1232} and Corollary \ref{cor:E-K-12-sqr},
\begin{align*}
\Cov_{1232}(\alphabh[])
\sim h^2 \bigg(\int_{R} f_{\eps[]}(y - \xb[1]^T \alphab[]) f_{\Xb[1]}(\xb[1]) f_{\Xb[3]}(\xb[3]) \ d\xb[1] d\xb[3] - \fY^2(y)\bigg),
\end{align*}
where the region of integration is $R =\{(\xb[1],\xb[3]) : (\xb[1] - \xb[3])^T \alphab[]=0\}$.
\end{proof}
\end{lemma}

%
%
%
%
\begin{lemma} \label{lem:Cov-1234}
It holds that
\begin{align*}
\Cov_{1234}(\alphabh[])
& \sim O(h^2 N^{-1}).
\end{align*}
\begin{proof}
Since $\Cov_{1234}(\alphabh[]) = E[K_{12}(\alphabh[])K_{34}(\alphabh[])] - E[K_{12}(\alphabh[])]^2$, by Lemma \ref{lem:E-K-1234} and Corollary \ref{cor:E-K-12-sqr},
\begin{align*}
\Cov_{1234}(\alphabh[])
& \sim h^2 N^{-1} \sigmaeps^2 \sum_{p_1,p_2=0}^J E[(\Phii[N]^{-1} \Xb[1])_{p_1}(\Phii[N]^{-1} \Xb[1])_{p_2}] \\
& \times E[f_{\eps[]}'(y - \Xb[1]^T \alphab[]) (\Xb[2] - \Xb[1])_{p_1}]E[f_{\eps[]}'(y - \Xb[3]^T \alphab[]) (\Xb[4] - \Xb[3])_{p_2}].
\end{align*}
\end{proof}
\end{lemma}

%
%
%
%
\begin{lemma} \label{lem:Cov-ijkl-1122}
It holds that
\begin{align*}
\Cov_{1122}(\alphabh[])
& = 0.
\end{align*}
\begin{proof}
By \eqref{eq:yt-simple},
\begin{align*}
\Cov_{1122}(\alphabh[])
&= E[K_{11}(\alphabh[])K_{22}(\alphabh[])] - E[K_{11}(\alphabh[])]^2 \\
&= E[K_h(y - \Xb[1]^T \alphab[] - \eps[1]) K_h(y - \Xb[2]^T \alphab[] - \eps[2])] - E[K_h(y - \Xb[1]^T \alphab[] - \eps[1])]^2 \\
&= 0.
\end{align*}
\end{proof}
\end{lemma}
The following lemmas relate to the covariances that turn out to be asymptotically negligible when the variance gets decomposed. 

%
%
%
%
\begin{lemma} \label{lem:Cov-ijkl-1231}
It holds that
\begin{align*}
\Cov_{ijkl}(\alphabh[])
& \sim O(h^2 N^{-1}),
\end{align*}
where $(i,j,k,l) \in \{(1,2,3,1),(1,2,2,1),(1,1,2,3),(2,1,1,3)\}$.
\end{lemma}

%
%
%
%
\begin{lemma} \label{lem:Cov-ijkl-1111}
It holds that
\begin{align*}
\Cov_{ijkl}(\alphabh[])
& \sim
\begin{cases}
O(h), & \quad \quad (i,j,k,l) = (1,1,1,1), \\
O(h^2), & \quad \quad (i,j,k,l) \in \{(1,1,1,2),(1,1,2,1)\}.
\end{cases}
\end{align*}
\end{lemma}
%
%
%
%
Now we are in a position to derive the asymptotic variance.
\begin{theorem}[] \label{thm:var-asy}
Under assumptions (A), (B), (C), (D), (E), and (F), the asymptotic variance of the multiple regression-enhanced convolution estimator is
\begin{align} \label{eq:variance}
\Var[\fhY(y)]
& = O\bigg(\frac{1}{hNL}\bigg) + O\bigg(\frac{1}{L}\bigg) + O\bigg(\frac{1}{N}\bigg).
\end{align}
\begin{proof}
The variance of the convolution estimator \eqref{eq:fhY-rewritten} can be decomposed as
\begin{align} \label{eq:var-proof-start}
\begin{split}
& \Var[\fhY(y)] \\
& = \frac{1}{(hNL)^2} \sum_{i=1}^L \sum_{j=1}^N \sum_{k=1}^L \sum_{l=1}^N \Cov_{ijkl}(\alphabh[]) \\
& = \frac{1}{(hNL)^2} \bigg(\sum_{i=1}^N \sum_{j=1}^N \sum_{k=1}^L \sum_{l=1}^N + \sum_{i={N+1}}^L \sum_{j=1}^N \sum_{k=1}^L \sum_{l=1}^N\bigg) \Cov_{ijkl}(\alphabh[]) \\
& = \frac{1}{(hNL)^2} \bigg(\sum_{i=1}^N \sum_{j=1}^N \sum_{k=1}^N \sum_{l=1}^N + \sum_{i=1}^N \sum_{j=1}^N \sum_{k=N+1}^L \sum_{l=1}^N + \sum_{i={N+1}}^L \sum_{j=1}^N \sum_{k=1}^N \sum_{l=1}^N + \sum_{i={N+1}}^L \sum_{j=1}^N \sum_{k=N+1}^L \sum_{l=1}^N\bigg) \Cov_{ijkl}(\alphabh[]) \\
& = \frac{1}{(hNL)^2} \bigg(\sum_{i=1}^N \sum_{j=1}^N \sum_{k=1}^N \sum_{l=1}^N + 2 \sum_{i={N+1}}^L \sum_{j=1}^N \sum_{k=1}^N \sum_{l=1}^N + \sum_{i={N+1}}^L \sum_{j=1}^N \sum_{k=N+1}^L \sum_{l=1}^N\bigg) \Cov_{ijkl}(\alphabh[]).
\end{split}
\end{align}
These sets of summations can be expressed as follows:
\begin{align} \label{eq:var-decomp-sum-A-B-C}
\begin{split}
\sum_{i=1}^N \sum_{j=1}^N \sum_{k=1}^N \sum_{l=1}^N \Cov_{ijkl}(\alphabh[]) \ & = \sum_{i=1}^{12} c_i^{(A)} \Psi_i^{(A)}, \\
\sum_{i={N+1}}^L \sum_{j=1}^N \sum_{k=1}^N \sum_{l=1}^N \Cov_{ijkl}(\alphabh[]) \ & = \sum_{i=1}^{5} c_i^{(B)} \Psi_i^{(B)}, \\
\sum_{i={N+1}}^L \sum_{j=1}^N \sum_{k=N+1}^L \sum_{l=1}^N \Cov_{ijkl}(\alphabh[]) \ & = \sum_{i=1}^{5} c_i^{(C)} \Psi_i^{(C)},
\end{split}
\end{align}
with
\begin{alignat*}{4}
&
c_{1}^{(A)} \sim N, & \quad \Psi_{1}^{(A)} = \Cov_{1111}(\alphabh[]), & \quad \quad
c_{1}^{(B)} \sim MN, & \quad \Psi_{1}^{(B)} = \Cov_{1222}(\alphabh[]), \\
&
c_{2}^{(A)} \sim N^2, & \quad \Psi_{2}^{(A)} = \Cov_{1212}(\alphabh[]), & \quad \quad
c_{2}^{(B)} \sim MN^2, & \quad \Psi_{2}^{(B)} = \Cov_{1223}(\alphabh[]), \\
&
c_{3}^{(A)} \sim N^3, & \quad \Psi_{3}^{(A)} = \Cov_{1213}(\alphabh[]), & \quad \quad
c_{3}^{(B)} \sim MN^3, & \quad \Psi_{3}^{(B)} = \Cov_{1234}(\alphabh[]), \\
&
c_{4}^{(A)} \sim 2N^3, & \quad \Psi_{4}^{(A)} = \Cov_{1231}(\alphabh[]), & \quad \quad
c_{4}^{(B)} \sim MN^2, & \quad \Psi_{4}^{(B)} = \Cov_{1232}(\alphabh[]), \\
&
c_{5}^{(A)} \sim N^3, & \quad \Psi_{5}^{(A)} = \Cov_{1232}(\alphabh[]), & \quad \quad
c_{5}^{(B)} \sim MN^2, & \quad \Psi_{5}^{(B)} = \Cov_{1233}(\alphabh[]), \\
& c_{6}^{(A)} \sim N^2, & \quad \Psi_{6}^{(A)} = \Cov_{1221}(\alphabh[]), & \quad \quad
c_{1}^{(C)} \sim MN, & \quad \Psi_{1}^{(C)} = \Cov_{1212}(\alphabh[]), \\
& c_{7}^{(A)} \sim 2N^2, & \quad \Psi_{7}^{(A)} = \Cov_{1112}(\alphabh[]), & \quad \quad
c_{2}^{(C)} \sim MN^2, & \quad \Psi_{2}^{(C)} = \Cov_{1213}(\alphabh[]), \\
& c_{8}^{(A)} \sim 2N^2, & \quad \Psi_{8}^{(A)} = \Cov_{1121}(\alphabh[]), & \quad \quad
c_{3}^{(C)} \sim M^2N, & \quad \Psi_{3}^{(C)} = \Cov_{1232}(\alphabh[]), \\
& c_{9}^{(A)} \sim N^4, & \quad \Psi_{9}^{(A)} = \Cov_{1234}(\alphabh[]), & \quad \quad
c_{4}^{(C)} \sim M^2N^2, & \quad \Psi_{4}^{(C)} = \Cov_{1234}(\alphabh[]), \\
& c_{10}^{(A)} \sim N^2, & \quad \Psi_{10}^{(A)} = \Cov_{1122}(\alphabh[]), \\
& c_{11}^{(A)} \sim N^3, & \quad \Psi_{11}^{(A)} = \Cov_{1123}(\alphabh[]), \\
& c_{12}^{(A)} \sim N^3, & \quad \Psi_{12}^{(A)} = \Cov_{2113}(\alphabh[]),
\end{alignat*}
It is worth noting that this is somewhat of a generalization of similar decompositions in \cite[Supp. Material]{stove2012convolution}. For example, the set $\{c_{i}^{(A)} \Psi_{i}^{(A)}\}_{i=1}^{12}$ can be related to the set $\{(S_i)\}_{i=1}^{12}$ in that work; see also \cite{saavedra2000estimation}. It suffices to evaluate the terms $\{\Psi_{i}^{(A)}\}_{i=1}^{12}$ since these terms correspond to equivalent terms in the sets $\{\Psi_{i}^{(B)}\}_{i=1}^{5}$ and $\{\Psi_{i}^{(C)}\}_{i=1}^{4}$. To be specific,
\begin{alignat*}{5}
& \Psi_{2}^{(A)} = \Psi_{1}^{(C)}, \quad \quad & \Psi_{3}^{(A)} = \Psi_{2}^{(C)}, \quad \quad & \Psi_{5}^{(A)} = \Psi_{4}^{(B)} = \Psi_{3}^{(C)}, \quad \quad & \Psi_{8}^{(A)} = \Psi_{1}^{(B)}, \\
& \Psi_{9}^{(A)} = \Psi_{3}^{(B)} = \Psi_{4}^{(C)}, \quad \quad & \Psi_{11}^{(A)} = \Psi_{5}^{(B)}, \quad \quad & \Psi_{12}^{(A)} = \Psi_{2}^{(B)}.
\end{alignat*}
By accounting for these correspondences, and using \eqref{eq:var-proof-start} and \eqref{eq:var-decomp-sum-A-B-C}, we can write the variance in the following form:
\begin{align} \label{eq:var-sum-d-Psi}
\Var[\fhY(y)]
= \sum_{i=1}^{12} d_{i} \Psi_{i}^{(A)},
\end{align}
where at leading-order,
\begin{alignat*}{8}
d_{1} &= \frac{c_{1}^{(A)}}{(hNL)^{2}} \ &\sim& \ \frac{1}{h^2NL^2}, \quad \quad & d_{7} &= \frac{2c_{7}^{(A)} }{(hNL)^{2}} \ &\sim& \ \frac{2}{h^2L^2}, \\
d_{2} &= \frac{(c_{2}^{(A)} + c_{1}^{(C)})}{(hNL)^{2}} \ &\sim& \ \frac{1}{h^2NL}, \quad \quad & d_{8} &= \frac{2(c_{8}^{(A)} + c_{1}^{(B)})}{(hNL)^{2}} \ &\sim& \frac{2}{h^2 NL}, \\
d_{3} &= \frac{(c_{3}^{(A)} + c_{2}^{(C)})}{(hNL)^{2}} \ &\sim& \ \frac{1}{h^2L}, \quad \quad & d_{9} &= \frac{(c_{9}^{(A)} + 2c_{3}^{(B)} + c_{4}^{(C)})}{(hNL)^{2}}  \ &\sim& \ \frac{1}{h^2}, \\
d_{4} &= \frac{2 c_{4}^{(A)}}{(hNL)^{2}} \ &\sim& \ \frac{2N}{h^2L^2}, \quad \quad & d_{10} &= \frac{c_{10}^{(A)}}{(hNL)^{2}} \ &\sim& \ \frac{1}{h^2L^2}, \\
d_{5} &= \frac{(c_{5}^{(A)} + 2c_{4}^{(B)} + c_{3}^{(C)})}{(hNL)^{2}} \ &\sim& \ \frac{1}{h^2N}, \quad \quad & d_{11} &= \frac{(c_{11}^{(A)} + 2c_{5}^{(B)})}{(hNL)^{2}} &\sim& \ \frac{1}{h^2L}, \\
d_{6} &= \frac{c_{6}^{(A)}}{(hNL)^{2}} \ &\sim& \ \frac{1}{h^2L^2}, \quad \quad & d_{12} &= \frac{(c_{12}^{(A)} + 2c_{2}^{(B)})}{(hNL)^{2}} &\sim& \ \frac{1}{h^2L}.
\end{alignat*}
Then, combing these expressions with the results for $\{\PsiA[i]\}_{i=1}^{12}$ derived in Lemmas \ref{lem:Cov-1212}, \ref{lem:Cov-1213}, \ref{lem:Cov-1232}, \ref{lem:Cov-1234}, \ref{lem:Cov-ijkl-1122}, \ref{lem:Cov-ijkl-1231}, and \ref{lem:Cov-ijkl-1111}, we find that
\begin{alignat*}{8}
d_{1} \PsiA[1] &\sim \frac{1}{h^2NL^2} O(h) \ &=& \ O\bigg(\frac{1}{hNL^2}\bigg), \quad \quad & d_{7} \PsiA[7] &\sim \frac{2}{h^2 L^2}O(h^2) \ &=& \ O\bigg(\frac{1}{L^2}\bigg), \\
d_{2} \PsiA[2] &\sim \frac{1}{h^2NL} O(h) \ &=& \ O\bigg(\frac{1}{hNL}\bigg), \quad \quad & d_{8} \PsiA[8] &\sim \frac{2}{h^2 NL}O(h^2) \ &=& \ O\bigg(\frac{1}{NL}\bigg), \\
d_{3} \PsiA[3] &\sim \frac{1}{h^2L}O(h^2) \ &=& \ O\bigg(\frac{1}{L}\bigg), \quad \quad & d_{9} \PsiA[9] &\sim \frac{1}{h^2}O(h^2N^{-1}) \ &=& \ O\bigg(\frac{1}{N}\bigg), \\
d_{4} \PsiA[4] &\sim \frac{2 N}{h^2L^2}O(h^2N^{-1}) \ &=& \ O\bigg(\frac{1}{L^2}\bigg), \quad \quad & d_{10} \PsiA[10] &\sim \frac{1}{h^2L^2} \cdot 0 \ &=& \ 0, \\
d_{5} \PsiA[5] &\sim \frac{1}{h^2N}O(h^2) \ &=& \ O\bigg(\frac{1}{N}\bigg), \quad \quad & d_{11} \PsiA[11] &\sim \frac{1}{h^2L}O(h^2N^{-1}) \ &=& \ O\bigg(\frac{1}{NL}\bigg), \\
d_{6} \PsiA[6] &\sim \frac{1}{h^2L^2}O(h^2N^{-1}) \ &=& \ O\bigg(\frac{1}{NL^2}\bigg), \quad \quad & d_{12} \PsiA[12] &\sim \frac{1}{h^2L}O(h^2N^{-1}) \ &=& \ O\bigg(\frac{1}{NL}\bigg).
\end{alignat*}
where we used the fact that $E[(\Phii[N]^{-1} \Xb[1])_{p_1}(\Phii[N]^{-1} \Xb[1])_{p_2}] = O(1)$ as $N \to \infty$ by Lemma \ref{lem:E-prod-PhiiN-Xi} for $\PsiA[9]$. Finally, as $N\to \infty$ and $h \to 0$, four of these terms are seen to dominate, that is to say,
\begin{align*}
\Var[\fhY(y)]
& \sim d_{2} \PsiA[2] + d_{3} \PsiA[3] + d_{5} \PsiA[5] + d_{9} \PsiA[9] \\
& = O\bigg(\frac{1}{hNL}\bigg) + O\bigg(\frac{1}{L}\bigg) + O\bigg(\frac{1}{N}\bigg) + O\bigg(\frac{1}{N}\bigg).
\end{align*}
\end{proof}
\end{theorem}

%
%
%
%
\begin{corollary}[] \label{cor:asy-mse}
Under assumptions (A), (B), (C), (D), (E), and (F), the asymptotic MSE of the multiple regression-enhanced convolution estimator is
\begin{align} \label{eq:mse}
\MSE[\fhY(y)]
&= \Bias[\fhY(y)]^2 - \Var[\fhY(y)]
= O(h^4) + O\bigg(\frac{1}{hNL}\bigg) + O\bigg(\frac{1}{L}\bigg) + O\bigg(\frac{1}{N}\bigg).
\end{align}
\end{corollary}

%
%
%
%
\begin{corollary}[] \label{cor:asy-mise}
Under assumptions (A), (B), (C), (D), (E), and (F), the asymptotic mean integrated square error (MISE) of the multiple regression-enhanced convolution estimator is
\begin{align} \label{eq:mise}
\MISE[\fhY]
&= \int \Bias[\fhY(y)]^2 - \Var[\fhY(y)] dy
= O(h^4) + O\bigg(\frac{1}{hNL}\bigg) + O\bigg(\frac{1}{L}\bigg) + O\bigg(\frac{1}{N}\bigg).
\end{align}
\end{corollary}

Note that in the case when there is no additional covariate information, that is when $M = 0$ and thus $L = N$, the mean square error reduces to $\MSE[\fhY(y)] = O(h^4) + O(h^{-1}N^{-2}) + O(N^{-1})$, which recovers a result established by St{\o}ve and Tj{\o}stheim \cite[Eq. (23)]{stove2012convolution}, and Escancianoa and Jacho-Ch{\'a}vez \cite[Eq. (3.4)]{escanciano2012n}.

Since $L = M+N$ appears in the variance \eqref{eq:variance} but not in the bias \eqref{eq:bias}, the presence of the supplemental sample of $M$ covariate observations leads to a direct reduction in the variance, but not the bias. However, the additional sample does have an indirect effect on the bias, since the asymptotically optimal bandwidth $h$ depends on $M$, as we show in the next section. Thus, the presence of the $M$ covariate observations in the additional sample ultimately leads to a reduction in the MSE and MISE through both the bias and the variance.


The effect that the amount of additional auxiliary data supplied to the convolution estimator has on the variance is clear from \eqref{eq:variance}. In the absence of additional auxiliary data, that is, when $M=0$ and thus $L=N$, we get
\begin{align} \label{eq:Var-C-0}
\Var[\fhY(y)] 
& = O\bigg(\frac{1}{hN^2}\bigg) + O\bigg(\frac{1}{N}\bigg).
\end{align}
When the number of additional covariate observations is on the order of the number of complete case observations, that is, when $M = O(N)$, we get
\begin{align} \label{eq:Var-C-times}
\Var[\fhY(y)] 
& = O\bigg(\frac{1}{hN(N+M)}\bigg) + O\bigg(\frac{1}{(N+M)}\bigg) + O\bigg(\frac{1}{N}\bigg),
\end{align}
which is just a rewritten version of the general expression \eqref{eq:variance} that makes the dependence on $M$ explicit. Finally, when the number of additional auxiliary data observations is much larger than the number of complete case observations, that is, when $N$ is large and fixed while $M \to \infty$, the first two terms in \eqref{eq:variance} vanish and we are left with
\begin{align} \label{eq:Var-C-gg-1}
\Var[\fhY(y)] 
& = O\bigg(\frac{1}{N}\bigg).
\end{align}
So, the presence of additional auxiliary data guarantees a reduction in the asymptotic variance.
No matter how many additional covariate observations are incorporated into the convolution estimator, however, the variance can't be reduced beyond $O(N^{-1})$ because $d_{5} \PsiA[5]$ and $d_{9} \PsiA[9]$ in \eqref{eq:var-sum-d-Psi} do not depend on the additional sample of $M$ covariate observations. This $O(N^{-1})$ uncertainty arises because the $O(N^{-1})$ uncertainty present in the underlying OLS estimator, which we recall was defined \eqref{eq:OLS-estimator} with respect to $N$ complete case observations, ultimately propagates into uncertainty in the convolution estimator.

Since the variance can't be reduced beyond $O(N^{-1})$, a saturation phenomenon arises. Eventually, as more and more additional auxiliary data observations are supplied to the convolution estimator, the improvement in accuracy will become completely negligible and the variance will saturate at $\Var[\fhY(y)] \sim d_{5} \PsiA[5] + d_{9} \PsiA[9]$.

%
%
%
%
\section{Bandwidth Selection and Convergence of the MSE with respect to $N$ and $M$} \label{sec:bandwidth-selection}

%
%
%
%
\begin{lemma}[] \label{lem:hopt}
The asymptotically optimal bandwidth $\hopt$ for the multiple regression-enhanced convolution estimator is
\begin{align} \label{eq:hopt}
\hopt = \bigg(\frac{\sigma_{K}}{\mu_K^2 \int (\fY''(y))^2 dy} \bigg)^{1/5} \frac{1}{(NL)^{1/5}}.
\end{align}
\begin{proof}
The MISE \eqref{eq:mise} depends on $h$ through the $h^2 \mu_K \fY''(y)/2$ term in the bias and the $d_2 \Psi_2^{(A)}$ term in the variance. Thus, the bandwidth $h$ which minimizes the MISE solves the following equation.
\begin{align*} 
0
&= \frac{\partial}{\partial h}\bigg(h^4 \frac{\mu_K^2}{4} \int (\fY''(y))^2 dy + \frac{\sigma_{K}}{hNL} \bigg)
= h^3 \mu_K^2 \int (\fY''(y))^2 dy - \frac{\sigma_{K}}{h^2NL}.
\end{align*}
\end{proof}
\end{lemma}

%
%
%
%
\begin{corollary}[] \label{cor:hopt-RP}
The asymptotically optimal bandwidth $\hopt$ for the multiple regression-enhanced convolution estimator can be written as
\begin{align} \label{eq:hopt-RP}
\hopt = \hoptRP L^{-1/5},
\end{align}
where
\begin{align*}
\hoptRP
&= \bigg(\frac{\sigma_{K}}{\mu_K^2 \int (\fY''(y))^2 dy} \bigg)^{1/5} \frac{1}{N^{1/5}},
\end{align*}
is the asymptotically optimal bandwidth for the Rosenblatt–Parzen density estimator \eqref{eq:kde} \cite[Eq. 3.21]{silverman1986density}.
\end{corollary}

%
%

Since $\fY''(y)$ is unknown, the optimal bandwidth formulas \eqref{eq:hopt} and \eqref{eq:hopt-RP} are not directly applicable. However, there are numerous methods in the literature for estimating $\hoptRP$, such as cross-validation \cite{rudemo1982empirical}, Silverman's rule of thumb \cite{silverman1986density}, and the plug-in approach of Sheather and Jones \cite{sheather1991reliable}. In any, case due to Corollary \ref{cor:hopt-RP}, we can employ an established means of choosing $\hoptRP$ and then scale it by $L^{-1/5}$ to obtain an estimate of the optimal bandwidth $\hopt$ for the convolution estimator. 

In the case when there is no additional covariate information, that is when $M = 0$ and thus $L = N$, convolution estimators already allow for reduced bias in comparison to the classical kernel density estimator \eqref{eq:kde}. For convolution estimators, $\hopt = O(N^{-2/5})$ which implies that $\Bias[\fhY(y)] = O(\hopt^2) + O(N^{-1}) =  O(N^{-4/5}) + O(N^{-1}) = O(N^{-4/5})$. For the classical kernel density estimator, on the other hand, $\hoptRP = O(N^{-1/5})$ which implies that $\Bias[\fhY(y)] = O((\hoptRP)^2) \sim O(N^{-2/5})$. The incorporation of additional covariate observations into the convolution estimator allows for an even greater reduction in bias, since the $h^2$ term in the bias in this case is $O((NL)^{-2/5})$ which is smaller than the $O(N^{-4/5})$ that arises in the usual case of the convolution estimator with no additional covariate information.

With the asymptotically optimal bandwidth in hand, we are now in a position to quantify the reduction in the MSE as the number of complete case observations $N$ increases, and the number of additional covariate observations $M$ increases.

%
%
%
%
\begin{lemma} \label{lem:MSE-N-infty}
For a fixed $M$, at the asymptotically optimal bandwidth, the MSE of the multiple regression-enhanced convolution estimator decays as $O(N^{-1})$ as $N \to \infty$.
\begin{proof}
As $N \to \infty$, we have that $\hopt = (NL)^{-1/5} \sim N^{-2/5}$, and thus
\begin{align*}
\MSE[\fhY(y)]
& = O(\hopt^4) + O\bigg(\frac{1}{\hopt NL}\bigg) + O\bigg(\frac{1}{L}\bigg) + O\bigg(\frac{1}{N}\bigg) \\
& = O\bigg(\frac{1}{N^{8/5}}\bigg) + O\bigg(\frac{N^{2/5}}{N^2}\bigg) + O\bigg(\frac{1}{N}\bigg) \\
& = O\bigg(\frac{1}{N^{8/5}}\bigg) + O\bigg(\frac{1}{N^{8/5}}\bigg) + O\bigg(\frac{1}{N}\bigg) \\
& = O\bigg(\frac{1}{N}\bigg).
\end{align*}
\end{proof}
\end{lemma}

%
%
%
%
\begin{lemma} \label{lem:MSE-M-infty}
For a large fixed $N$, at the asymptotically optimal bandwidth, the MSE of the multiple regression-enhanced convolution estimator decays as $O(M^{-4/5})$ towards an $O(N^{-1})$ constant as $M \to \infty$.
\begin{proof}
For $N$ large, we have that
\begin{align*}
\MSE[\fhY(y)]
& = O(\hopt^4) + O\bigg(\frac{1}{\hopt NL}\bigg) + O\bigg(\frac{1}{L}\bigg) + O\bigg(\frac{1}{N}\bigg) \\
& = O\bigg(\frac{1}{(NL)^{4/5}}\bigg) + O\bigg(\frac{(NL)^{1/5}}{NL}\bigg) + O\bigg(\frac{1}{L}\bigg) + O\bigg(\frac{1}{N}\bigg) \\
& = O\bigg(\frac{1}{(NL)^{4/5}}\bigg) + O\bigg(\frac{1}{L}\bigg) + O\bigg(\frac{1}{N}\bigg).
\end{align*}
The $O(N^{-1})$ third term acts as the saturation threshold since it doesn't change with $M$. Now, with $N$ fixed, the first term on the last line decays as $O((NL)^{-4/5}) = O((N^2 + NM)^{-4/5}) = O(M^{-4/5})$. This term dominates the second term on the last line which decays as $O(L^{-1}) = O((N + M)^{-1}) = O(M^{-1})$.
\end{proof}
\end{lemma}
It is worth highlighting the fact that the convergence of the MSE of the multiple regression-enhanced convolution estimator with respect to $N$ and $M$ is independent of the dimensionality of the covariates. This is in contrast to the case of Nadaraya-Watson which suffers from the curse of dimensionality \cite[Sec. 6]{stove2012convolution}. Therefore, if the data is well-fit by a multiple regression model, possibly after some variable transformations, and the dimensionality of the covariate vector is significant, it is recommended to use the multiple regression-enhanced convolution estimator to achieve $O(N^{-1})$ convergence.

By characterizing the decay of the MSE with respect to $M$, Lemma \ref{lem:MSE-M-infty} resolves the question regarding precisely how much the accuracy of density estimates can be enhanced through the incorporation of additional covariate information into a convolution estimator. Since increasing the number of complete case observations causes the MSE to decay as $O(N^{-1})$ towards zero, while increasing number of additional covariate observation causes the MSE to decay as $O(M^{-4/5})$ towards an $O(N^{-1})$ constant, we see that supplying more covariate observations to the convolution estimator is not quite as effective as supplying more complete case observations. However, supplying additional covariate observations can still provide a very significant performance improvement, as will be demonstrated in the numerical simulations in Section \ref{sec:numerical-simulations}. This is good news because in many in practical applications it may be difficult if not downright impossible to obtain additional response observations, whereas it can often be very straightforward to obtain large amounts of additional covariate data. 
\section{Efficient Computational Implementation} \label{sec:numerical-implementation}
Since the expression for the multiple regression-enhanced convolution estimator \eqref{eq:fhY} involves summing over the entire set of $L$ covariate observations, and also the set of $N$ residuals from the multiple regression model \eqref{eq:samnple-mult-reg}, the evaluation of the convolution estimator can be very time consuming in comparison to the evaluation of the classical density estimator \eqref{eq:kde} which features only a single summation over the set of response observations. In particular, in the case of large datasets, or in applications involving cross-validation or bootstrapping, the computational costs can become prohibitive.

To reduce computational times, in this section we present a Fast Gauss Transform (FGT)-based acceleration algorithm. This algorithm utilizes the high-performance C++ library FIGTree \cite{morariu2008automatic} which combines the (Improved) Fast Gauss Transform \cite{greengard1991fast} with Approximate Nearest Neighbor searching \cite{arya1993approximate} to efficiently evaluate the Gauss transform. The Gauss transform $G$ is defined as
\begin{align} \label{eq:FGT}
G_h(y_v,\{x_i\}_{i=1}^L,\{q_i\}_{i=1}^L) = \sum_{i=1}^L q_i e^{-((y_v-x_i)/h)^2},
\end{align}
where $\{q_j\}_{i=1}^L$ is a set of coefficients, $\{y_v\}_{v=1}^V$ is a set of target points, and $\{x_i\}_{i=1}^L$ is a set of source points.
The computational complexity involved in directly evaluating this expression at the $V$ target points is $O(VL)$. The FGT reduces the computational complexity to $O(V+L)$.

We need to evaluate the convolution estimator \eqref{eq:fhY} at the set of $V$ target points:
\begin{align} \label{eq:fhY-target-points}
\fhY(y_v) = \frac{1}{hNL} \sum_{i=1}^L \sum_{j=1}^N K_h(y_v - \Xb[i]^T \alphabh[]-\epsh[j]).
\end{align}
Since this expression features an extra summation compared to \eqref{eq:FGT}, a naive computational implementation results in a complexity of $O(VLN)$. By employing the FGT, the complexity can be reduced to $O(VN + L)$. With the complexity now scaling linearly with the number of complete case observations and additional covariate observations, the algorithm proposed in this section is particularly effective for accelerating evaluations when the number of additional covariate observations is potentially orders of magnitude larger than the number of complete cases observations, such as in the case of a difficult to measure response variable.

FIGTree cannot be used directly for the computation of \eqref{eq:fhY-target-points}, since FIGTree accelerates the evaluation of the single summation in \eqref{eq:FGT}. Therefore, we need to transform \eqref{eq:fhY} into a single summation expression. This can be achieved by stacking the set of $N$ residuals and $V$ target points into a single set of $VN$ artificial target points. Define the artificial target points $\{z_k\}_{k=1}^{VN}$ by
\begin{align} \label{eq:z-k}
z_k = y_{\floor{(p-1)/N}+1} - \hat{\varepsilon}_{((p-1 \bmod N)+1}, \quad \quad p = 1, \dots, VN.
\end{align}
Then \eqref{eq:fhY-target-points} can be rewritten as
\begin{align} \label{eq:fhY-target-points-figtree}
\fhY(y_v) = \frac{1}{\sqrt{2\pi}}\frac{1}{hNL} \sum_{p=(v-1)N+1}^{vN} G_h(z_p,\{\Xb[i]^T \alphabh[]\}_{i=1}^L,\{1\}_{i=1}^L).
\end{align}
Now, $\{G_h(z_k,\{\Xb[i]^T \alphabh[]\}_{i=1}^L,\{1\}_{i=1}^L)\}_{k=1}^{VN}$ can be evaluated using FIGTree. Once the evaluations have been performed at the artificial target points, the convolution estimator evaluations at the actual target points $\{y_v\}_{v=1}^V$ can be recovered using \eqref{eq:fhY-target-points-figtree}; see Algorithm \ref{alg:cnv-est}.

In Table \ref{tab:comp-times}, we present computational times for the evaluation of the multiple regression-enhanced convolution estimator \eqref{eq:fhY}, using a variety of approaches in the case of $N=100$ complete case observations and $V = 50$ target points, as a progressively larger number of additional covariate observations is supplied to the estimator. The evaluation approaches are as follows.
\begin{enumerate}[(i)]
\item Naive (R): Evaluating the density estimator using for loops in R.
\item Naive (C++): Evaluating the density estimator using for loops in C++.
\item FFT: Stacking the $L$  variables in the set $\{\Xb[i]^T \alphab[]\}_{i=1}^L$ and the $N$ variables in the set $\{\epsh[j]\}_{i=1}^N$ in \eqref{eq:fhY-target-points} into a single set of length $LN$, and then evaluating the density estimator using the \texttt{density()} function in R, which employs the Fast Fourier Transform (FFT) to accelerate computations \cite{deng2011density}.
\item FGT: The FGT acceleration technique presented above.
\end{enumerate}
Naive evaluation in R is very slow, which is to be expected since R is an interpreted language. For $M=12800$ additional covariate observations, naive evaluation in C++ is about $33$ times faster than naive evaluation in R. Performing the evaluations using the FFT-accelerated \texttt{density()} function in R is about $1.8$ times faster than performing the evaluations using Naive (C++).

At $M=100$, the FGT acceleration technique introduced above is approximately $159, 5.9$, and $4.6$ times faster than performing the evaluations using Naive (R), Naive (C++), and the FFT, respectively. As $M$ increases, the acceleration becomes even more pronounced; at $M=12800$, the FGT acceleration technique is approximately $2722, 81$, and $41$ times faster than performing the evaluations using Naive (R), Naive (C++), and the FFT, respectively. At this value of $M$, the FGT evaluation takes about 22 milliseconds whereas the FFT evaluation takes almost $1$ second.

Note that this algorithm can easily be adapted to other types of convolution estimators such as the Nadaraya-Watson-based convolution estimator; we can simply replace the set $\{\Xb[i]^T \alphabh[]\}_{i=1}^L$ with the analogous set $\{\hat{m}(\Xb[i])\}_{i=1}^L$, where $\hat{m}$ is the Nadaraya-Watson estimator of the nonlinear regression function $m$ in \eqref{eq:reg-nonlin}.
\renewcommand{\arraystretch}{1.2}
\begin{table}[htb!]
\footnotesize
\centering
\begin{tabular}[t]{|l|rrrrrrrrr|}
\hline
M        & 0 & 100 & 200 & 400 & 800 & 1600 & 3200 & 6400 & 12800 \\
\hline
Naive (R) & 4.86e-01 & 1.10e+00 & 1.39e+00 & 2.30e+00 & 4.23e+00 & 7.72e+00 & 1.53e+01 & 2.96e+01 & 6.01e+01 \\
Naive (C++) & 1.97e-02 & 4.12e-02 & 5.65e-02 & 9.21e-02 & 1.60e-01 & 2.87e-01 & 5.13e-01 & 9.62e-01 & 1.78e+00 \\
FFT   & 6.98e-03 & 3.18e-02 & 1.87e-02 & 3.65e-02 & 6.27e-02 & 1.16e-01 & 2.17e-01 & 4.31e-01 & 8.99e-01 \\
FGT   & 4.99e-03 & 6.98e-03 & 6.02e-03 & 6.98e-03 & 6.03e-03 & 7.92e-03 & 9.99e-03 & 1.49e-02 & 2.21e-02 \\
\hline
\end{tabular}
\caption{\label{tab:comp-times} Computational times (seconds) for evaluations of the multiple regression-enhanced convolution estimator as a progressively larger number of additional covariate observations are supplied to the estimator, using a variety of evaluation approaches. The FGT-based algorithm is much faster than the other approaches in all cases, with the reduction in computational times becoming even more pronounced as $M$ increases.}
\end{table}

\begin{algorithm}[htb]
\caption{Multiple regression-enhanced convolution estimator \label{alg:cnv-est}}
\begin{enumerate}
\item Fit a regression model to the complete case dataset $\{(\Y[i], \Xb[i])\}_{i=1}^N$ to obtain the OLS estimator $\alphabh[]$ and the residuals $\{\epsh[j]\}_{i=1}^N$.
\item Combine the $N$ covariate observations in the complete case dataset and the $M$ covariate observations in the additional dataset into a single dataset $\{\Xb[i]\}_{i=1}^L$.
\item Generate the artificial target points $\{z_k\}_{k=1}^{VN}$ using \eqref{eq:z-k}.
\item Compute $\{G_h(z_k,\{\Xb[i]^T \alphabh[]\}_{i=1}^L,\{1\}_{i=1}^L)\}_{k=1}^{VN}$ by supplying the covariate observations $\{\Xb[i]\}_{i=1}^L$ and the artificial target points $\{z_k\}_{k=1}^{VN}$ to FIGTree which evaluates \eqref{eq:FGT}.
\item Convert $\{G_h(z_k,\{\Xb[i]^T \alphabh[]\}_{i=1}^L,\{1\}_{i=1}^L)\}_{k=1}^{VN}$ into density estimate evaluations at the actual target points $\{\fhY(y_v)\}_{v=1}^V$ using \eqref{eq:fhY-target-points-figtree}.
\end{enumerate}
\end{algorithm}

\section{Numerical Simulations} \label{sec:numerical-simulations}
In this section we present some numerical simulations to investigate the accuracy of the multiple regression-enhanced convolution estimator. Let $\{y_v\}_{v=1}^V$ be a uniformly spaced set of target points, where $V = 128$. To compute the MISE of Rozenblatt-Parzen estimator and the multiple regression-enhanced convolution estimator we use a reference solution $\fY$ which is obtained by estimating the density of the response variable with the Rozenblatt-Parzen estimator using a very large sample size, $N = 10^6$.

Denote by $\fhY^{(p)}(y_v)$ a realization of an estimate given by a density estimator at the target point $y_v$.  We approximate the integrated square error (ISE) of the realization $\fhY^{(p)}$ by the Riemann sum
\begin{align*}
\ISE[\fhY^{(p)}]
&= (v_2-v_1)\sum_{v=2}^V (\fhY^{(p)}(y_v) - \fY(y_v))^2,
\end{align*}
where $P = 500$. The MISE of the estimator $\fhY$ is then approximated by
\begin{align*}
\MISE[\fhY]
&= \frac{1}{P} \sum_{p=1}^P \ISE[\fhY^{(p)}].
\end{align*}

Denote by $\fhY^{(RP)}$ the Rosenblatt–Parzen density estimator, and by $\fhY^{(MR)}$ the multiple regression-enhanced convolution estimator. Since the performance of convolution estimators with respect to the number of complete cases observations $N$ has already been analyzed in papers such as \cite{stove2012convolution,li2016n,escanciano2012n,muller2012estimating}, in this work we are more concerned with the performance with respect to the number of additional covariate observations $M$. In particular, we are interested in investigating by how much $\MISE[\fhY^{(MR)}]$ can be reduced as the number of additional covariate observations supplied to the convolution estimator increases. Denote by $\tau = M/N$ the ratio of additional covariate observations to complete case observations.

For the convolution estimator, we use bandwidths given by the asymptotically optimal bandwidth formula \eqref{eq:hopt-RP}, which involves scaling the corresponding asymptotically optimal bandwidths for the Rozenblatt-Parzen estimator. For the Rozenblatt-Parzen estimator itself, we use the Sheather-Jones method of bandwidth selection \cite{sheather1991reliable}.

\subsection{Single peaked negatively skewed distribution}

\begin{figure}[ht!]
	\centering
	\captionsetup[subfloat]{labelformat=empty}
	\begin{tabular}{ccc}
	\subfloat[(i) $\tau=0$]{\includegraphics[scale=0.6]{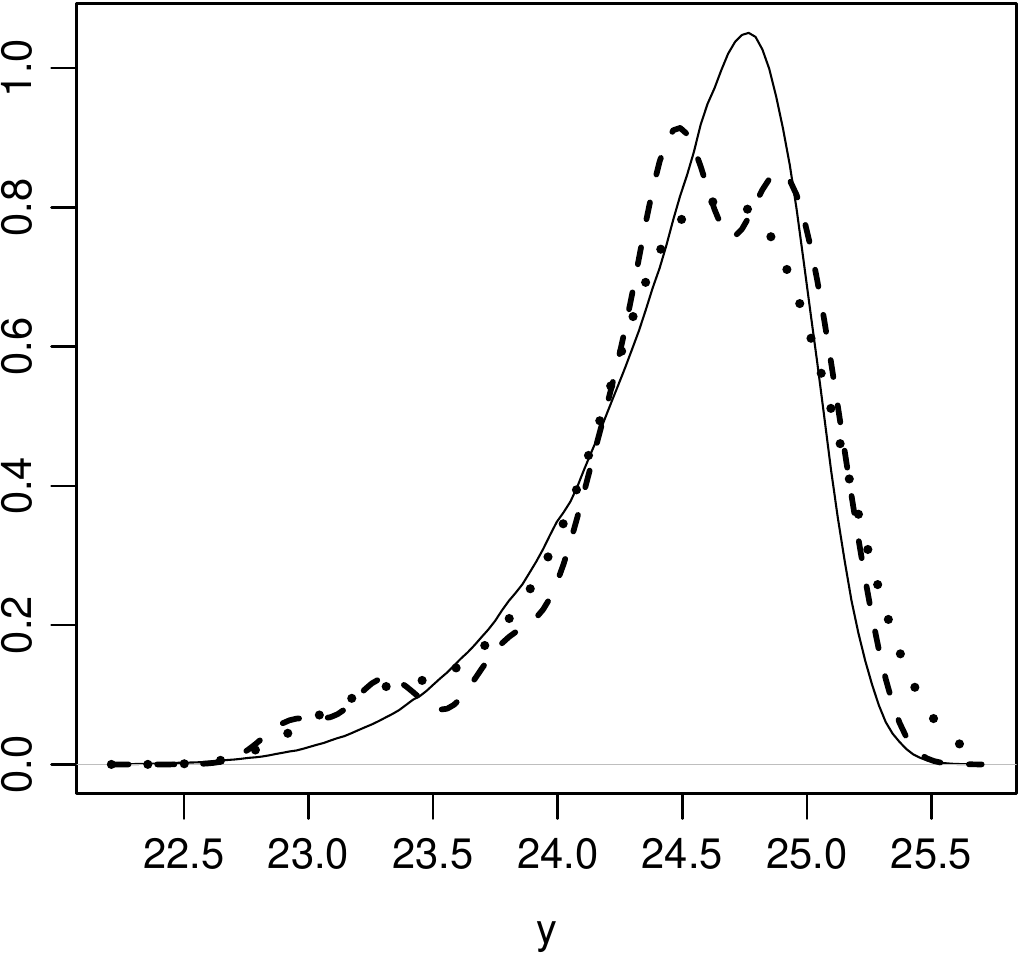}} \hspace{2em}
	\subfloat[(ii) $\tau=4$]{\includegraphics[scale=0.6]{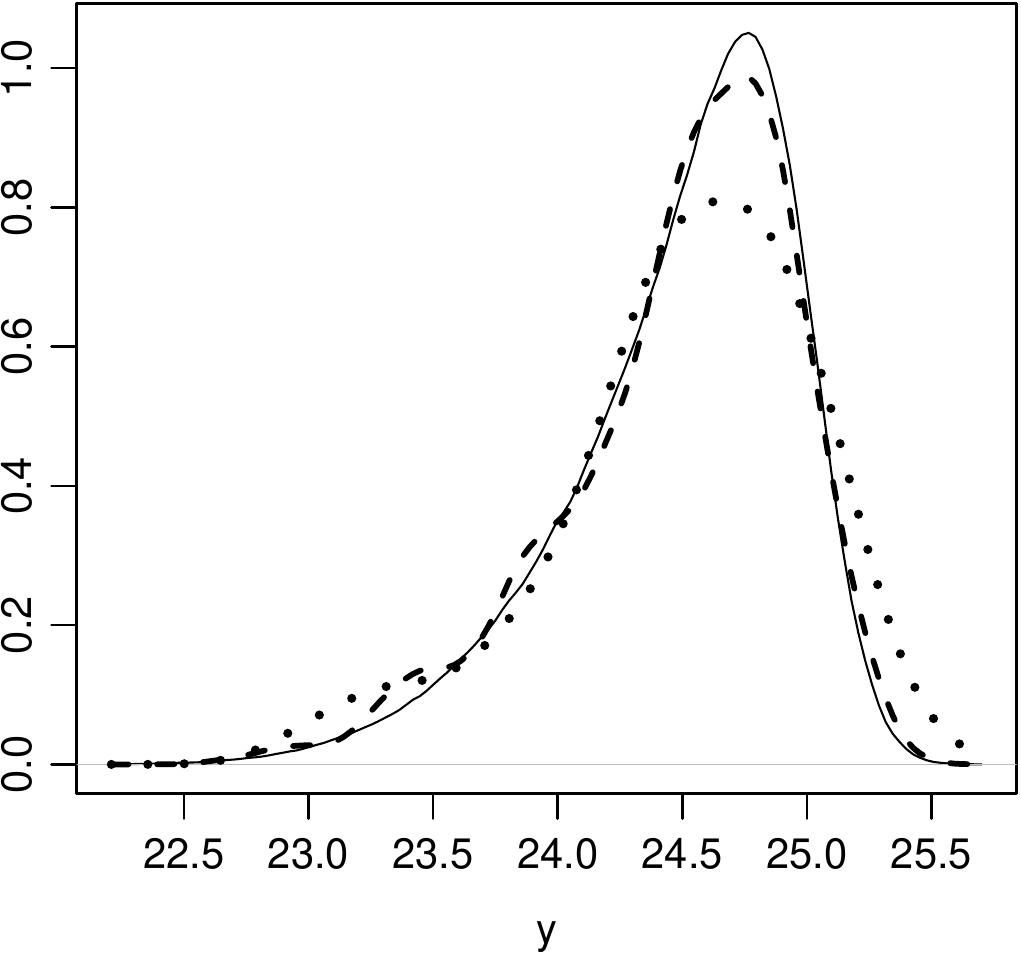}} \\
	\subfloat[(iii) $\tau=16$]{\includegraphics[scale=0.6]{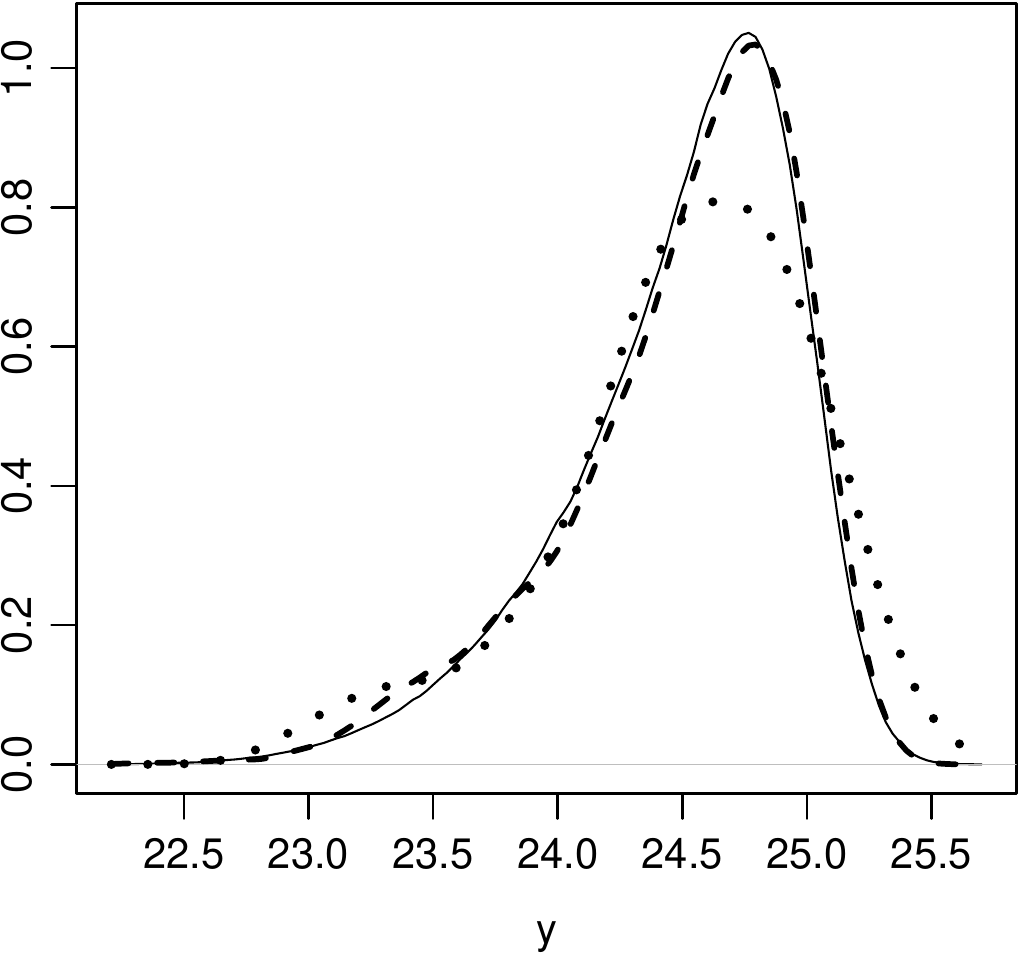}} \hspace{2em}
	\subfloat[(iv) $\tau=64$]{\includegraphics[scale=0.6]{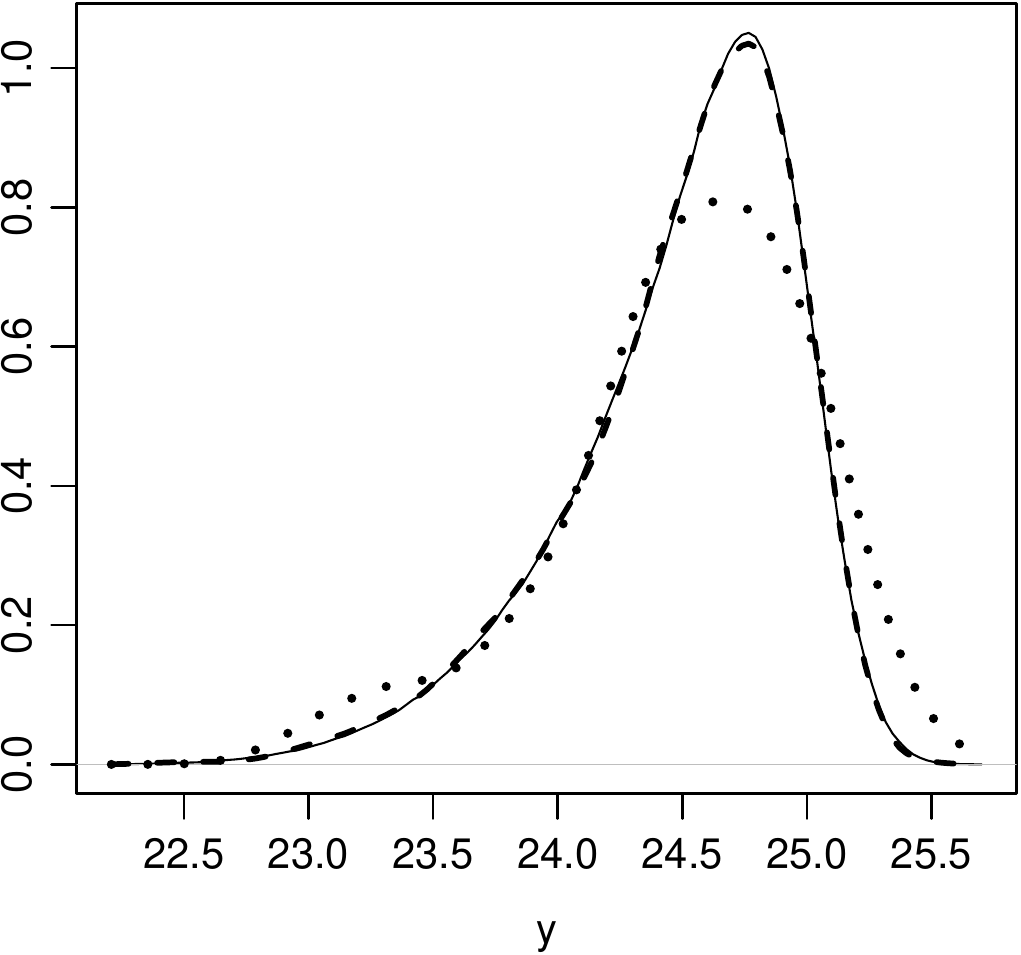}}
	\end{tabular}
	\caption{Typical realizations of estimates for the density of the response variable $Y$ in \eqref{eq:mult-reg-skewed} given by the Rosenblatt–Parzen density estimator $\fhY^{(RP)}$ (dotted line), and the multiple regression-enhanced convolution estimator $\fhY^{(MR)}$ (dashed line). The black solid line is the true density $\fY$. \label{fig:typ-real-skew}}
\end{figure}

Consider the regression model
\begin{align} \label{eq:mult-reg-skewed}
Y = \alpha_0 + \alpha_1 X_1 + \alpha_2 X_2 + \eps[],
\end{align}
where $(\alpha_0,\alpha_1,\alpha_2) = (1,3,3)$, with
\begin{align*}
X_1 \sim \beta(5,1), \quad \quad
X_2 \sim \mathcal{N}(7,0.05), \quad \quad
\eps[] & \sim \mathcal{N}(0,0.1).
\end{align*}
The density $\fY$ of $Y$ is in this case is single-peaked and negatively skewed. In Figure \ref{fig:typ-real-skew}, we plot $\fY$ along with typical realizations of density estimates given by $\fhY^{(RP)}$ and $\fhY^{(MR)}$, when the complete case dataset features $N=100$ observations, while the ratio of additional covariate observations supplied to the convolution estimator progressively increases, $\tau \in \{0,4,16,64\}$. In subplot (i), both $\fhY^{(RP)}$ and $\fhY^{(MR)}$ provide poor estimates, particularly in the under-smoothed tail region.
The estimate $\fhY^{(MR)}$ for $\tau=0$ is no better than $\fhY^{(RP)}$. If anything, its actually worse since it is bimodal when the true distribution is unimodal. However, as $\tau$ increases, $\fhY^{(MR)}$ approaches the true distribution; the problematic tail region gets smoothed out and $\fhY^{(MR)}$ becomes unimodal. At $\tau = 64$, $\fhY^{(MR)}$ provides a very accurate estimate of $\fY$. 

In Table \ref{tab:skew} we report MISE results for both $\fhY^{(RP)}$ and $\fhY^{(MR)}$. Three complete cases samples sizes are considered, $N \in \{50,100,200\}$. For each sample size, we compute the MISE when the ratio of additional covariate observations to complete case observations is $\tau \in \{0,2,4,8,16,32,64,128,256,512\}$. Note that for all complete case sample sizes, $\MISE[\fhY^{(MR)}] \approx \MISE[\fhY^{(RP)}]$ when $\tau = 0$. Now consider, for example, the $N=200$ case. By the time $\tau = 128$, $\MISE[\fhY^{(MR)}]$ is about $23$ times smaller than $\MISE[\fhY^{(RP)}]$, which demonstrates that a very substantial reduction in MISE is achievable through the incorporation of additional covariate observations. The change in the MISE as $\tau$ increase from $\tau = 128$ to $\tau = 512$ is negligible as saturation has occurred by this stage.

\begin{figure}[h]
	\centering
	\includegraphics[scale=0.85]{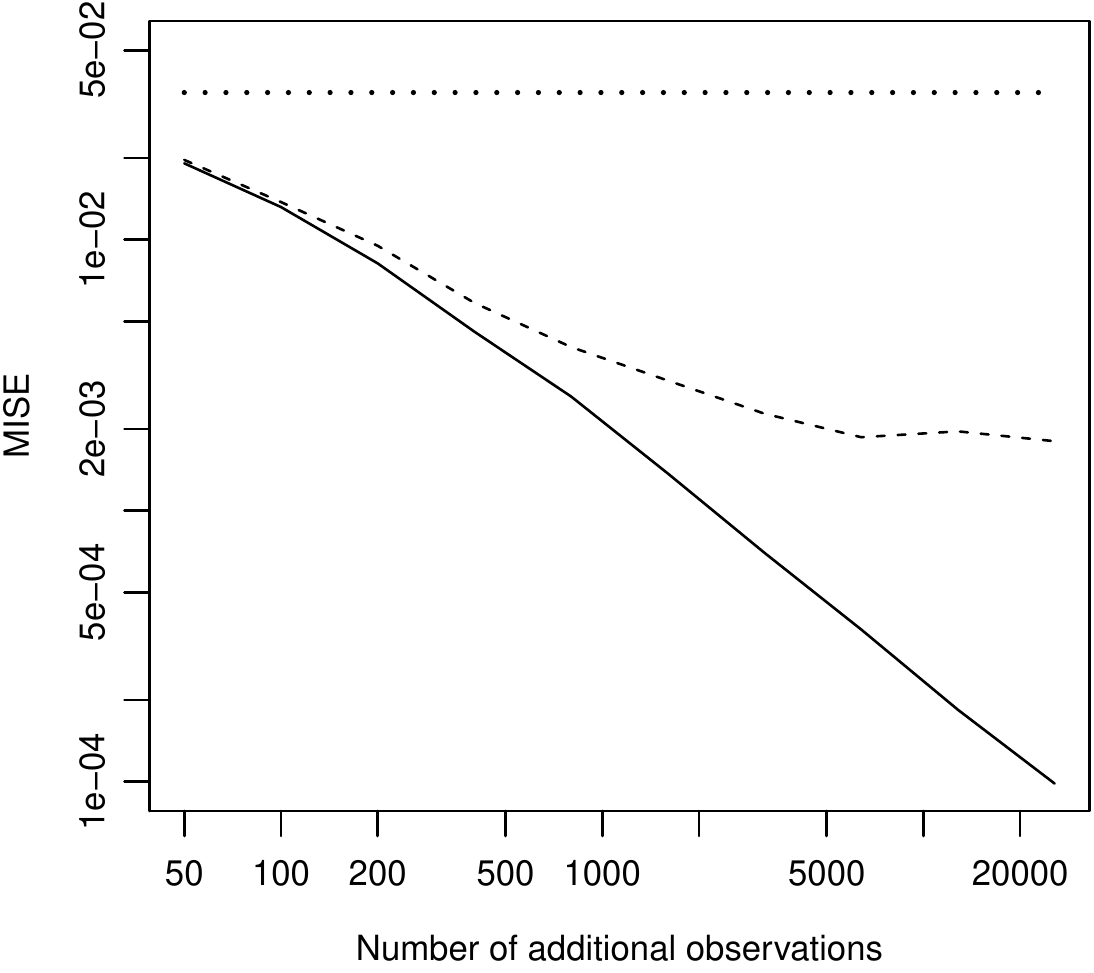}
	\caption{The convergence of the MISE as additional observations are supplied to the convolution estimator. The $\MISE[\fhY^{(MR)}]$ results from Table \ref{tab:skew} in which $N=50$ complete case observations are supplemented by an additional $M = \tau N$ covariate observations are represented by the dashed line. The solid black line is the MISE of $\fhY^{(MR)}$ when the convolution estimator is supplied with a progressively larger number of additional complete case observations. The corresponding $\MISE[\fhY^{(RP)}]$ result is the dotted line.
\label{fig:skewed-mise-cvg}}
\end{figure}

To understand the difference between incorporating additional covariate observations versus incorporating additional complete case observations in the context of MISE reduction, consider the loglog plot in Figure \ref{fig:skewed-mise-cvg}. The $\fhY^{(MR)}$ MISE results from Table \ref{tab:skew} for the case of $N=50$ are represented by the dashed line in this plot. The fact that $\MISE[\fhY^{(MR)}]$ is slightly higher at $\tau = 256$ compared to $\tau = 128$ and $\tau = 512$ is just a numerical artefact of the convergence flat-lining once the saturation threshold has been reached. As $M$ increases, this MISE is converging asymptotically as $O(M^{-4/5})$ towards an $O(N^{-1})$ constant. The solid black line is the MISE of $\fhY^{(MR)}$ when the initial dataset of $N=50$ complete case observations is supplemented with a progressively larger number of additional complete case observations as opposed to additional covariate observations. This line decays asymptotically as $O(N^{-1})$ towards zero. The corresponding $\fhY^{(RP)}$ MISE result is also shown for reference as a dotted line. Clearly, incorporating additional complete cases observations is more effective than incorporating additional covariate observations. Nevertheless, incorporating more covariate observations still allows for a substantial reduction in the MISE, which is very useful in situations in which obtaining more covariate observations is straightforward while obtaining more complete cases observations may be impossible.

\renewcommand{\arraystretch}{1.4}
\setlength{\tabcolsep}{.4em}
\begin{table}[ht]
\footnotesize
\centering
\begin{tabular}[t]{|l|c|cccccccccc|}
\hline
      & $\MISE[\fhY^{(RP)}]$ & \multicolumn{10}{c|}{$\MISE[\fhY^{(MR)}]$} \\
      &                & $\tau=0$ & $\tau=2$ & $\tau=4$ & $\tau=8$ & $\tau=16$ & $\tau=32$ & $\tau=64$ & $\tau=128$ & $\tau=256$ & $\tau=512$ \\
\hline
$N=50$  &  $3.49\text{e-}02$               & $3.57\text{e-}02$ & $1.38\text{e-}02$ & $9.50\text{e-}03$ & $5.85\text{e-}03$ & $4.02\text{e-}03$ &   $3.02\text{e-}03$     & $2.28\text{e-}03$ & $1.87\text{e-}03$ & $1.96\text{e-}03$ & $1.80\text{e-}03$ \\
$N=100$ &  $2.04\text{e-}02$               & $1.98\text{e-}02$ & $7.33\text{e-}03$ & $4.76\text{e-}03$ & $3.06\text{e-}03$ & $2.04\text{e-}03$ &   $1.44\text{e-}03$     & $1.21\text{e-}03$ & $9.01\text{e-}04$ & $9.27\text{e-}04$ & $8.71\text{e-}04$ \\
$N=200$ &  $1.16\text{e-}02$               & $9.61\text{e-}03$ & $3.82\text{e-}03$ & $2.44\text{e-}03$ & $1.62\text{e-}03$ & $1.01\text{e-}03$ &   $7.84\text{e-}04$     & $6.03\text{e-}04$ & $5.17\text{e-}04$ & $4.69\text{e-}04$ & $4.36\text{e-}04$ \\
\hline
\end{tabular}
\caption{\label{tab:skew} MISE results for the Rozenblatt-Parzen estimator and the multiple regression-enhanced convolution estimator for the density of the response variable in the regression model \eqref{eq:mult-reg-skewed}.}
\end{table}

\subsection{Response variable with multimodal distribution}
\begin{figure}[h]
	\centering
	\captionsetup[subfloat]{labelformat=empty}
	\begin{tabular}{ccc}
	\subfloat[(i) $\tau=0$]{\includegraphics[scale=0.6]{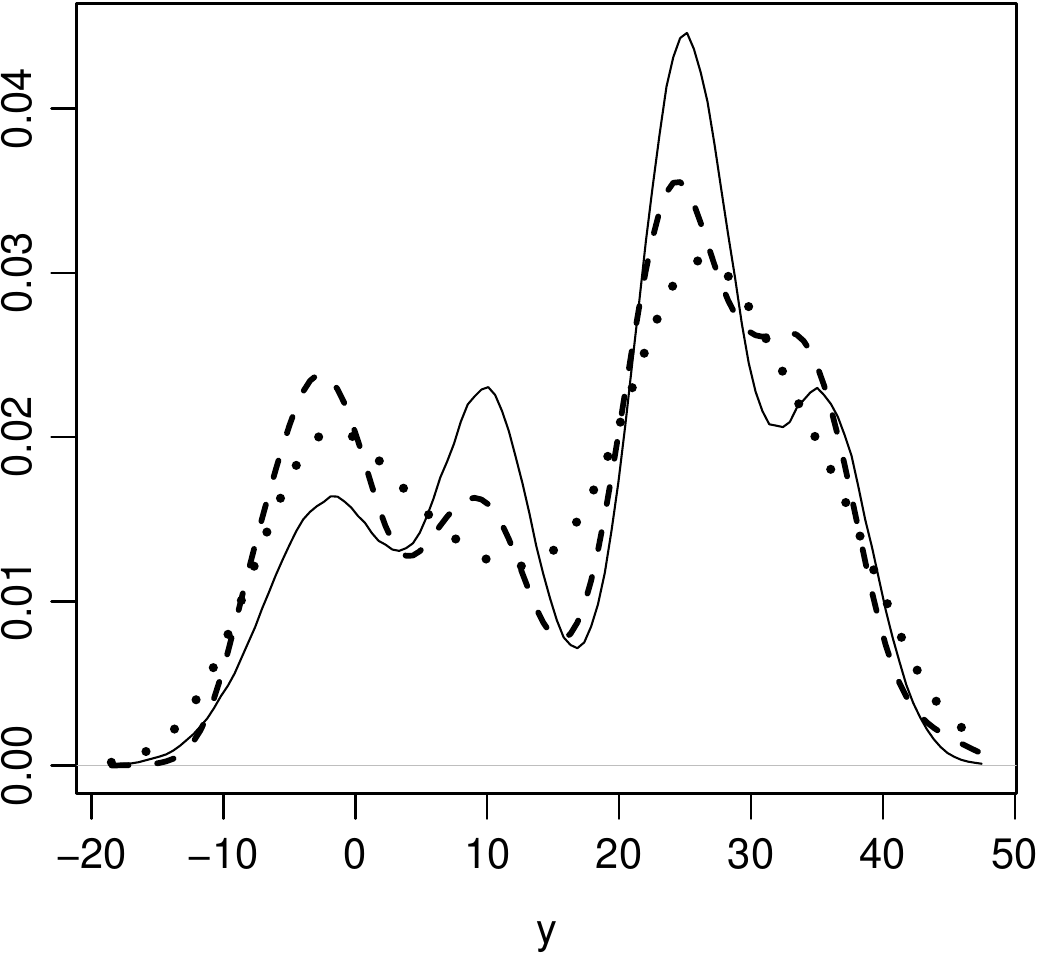}} \hspace{2em}
	\subfloat[(ii) $\tau=4$]{\includegraphics[scale=0.6]{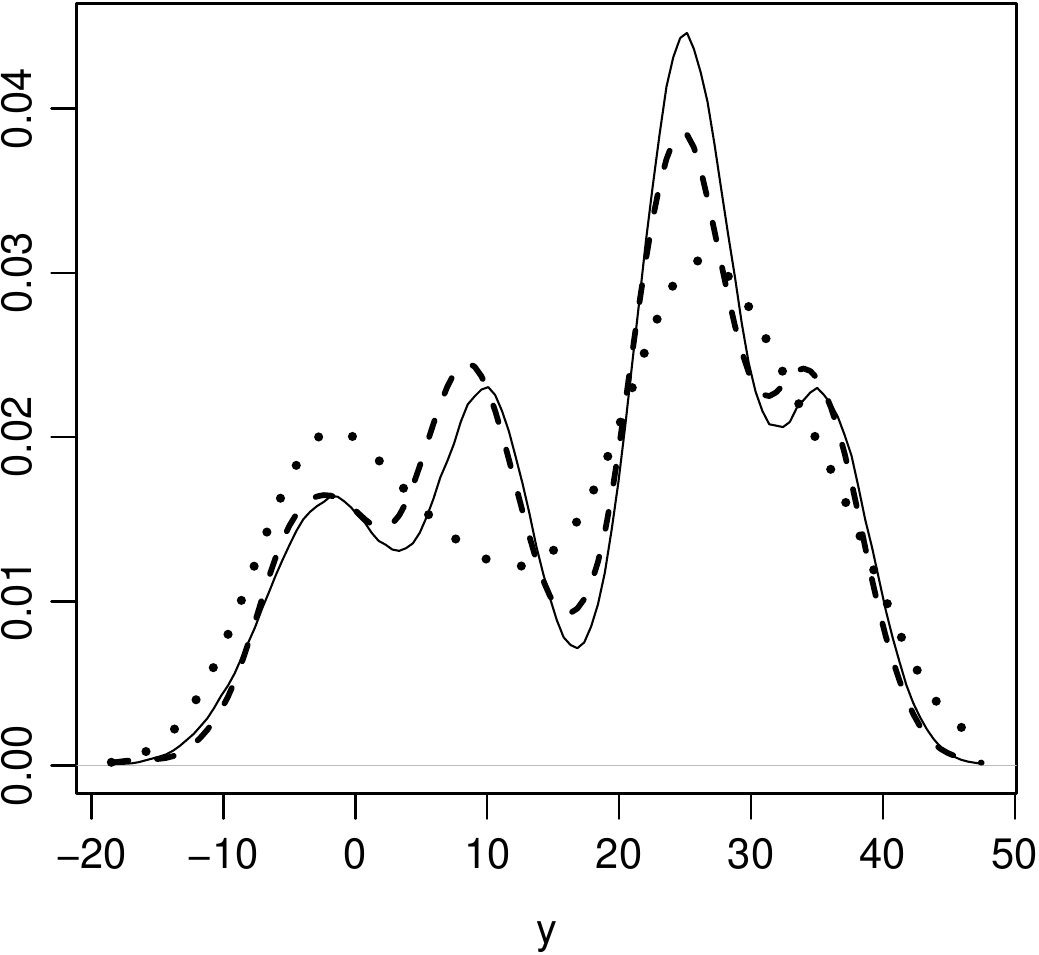}} \\
	\subfloat[(iii) $\tau=16$]{\includegraphics[scale=0.6]{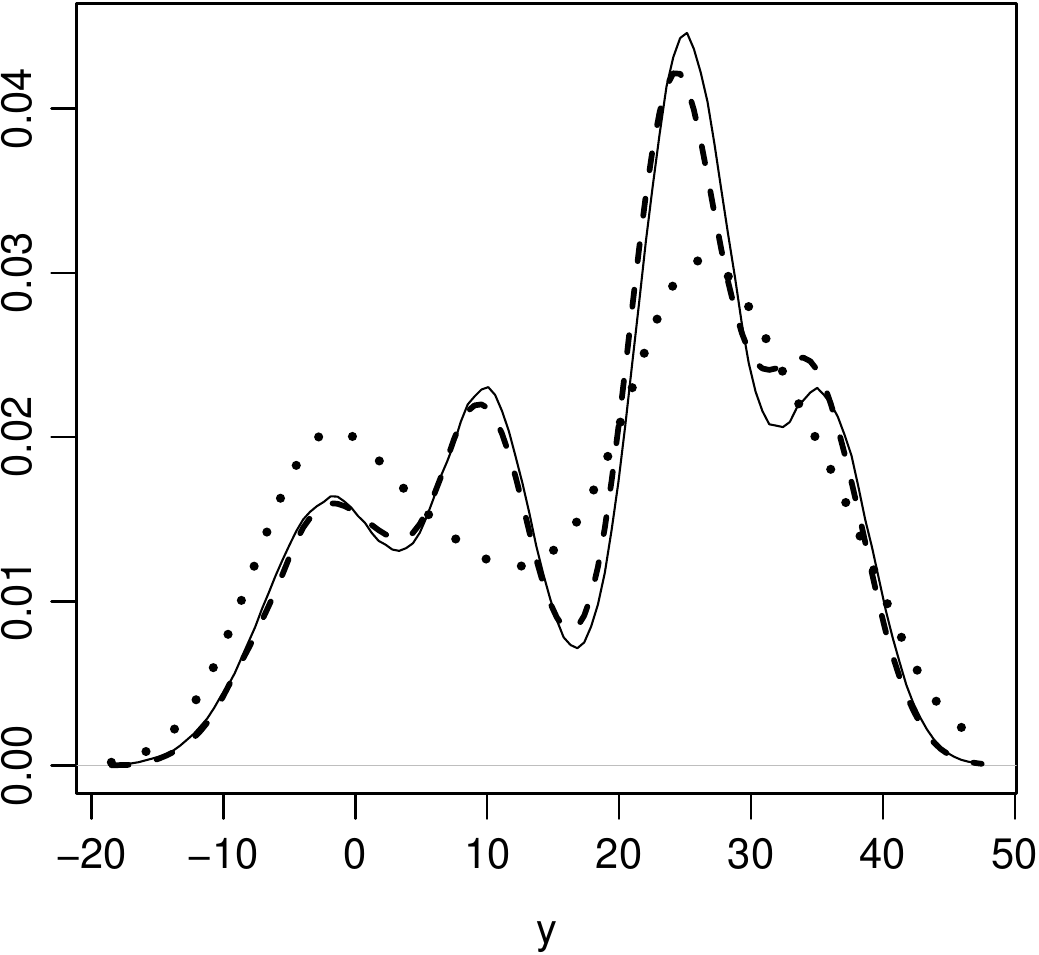}} \hspace{2em}
	\subfloat[(iv) $\tau=64$]{\includegraphics[scale=0.6]{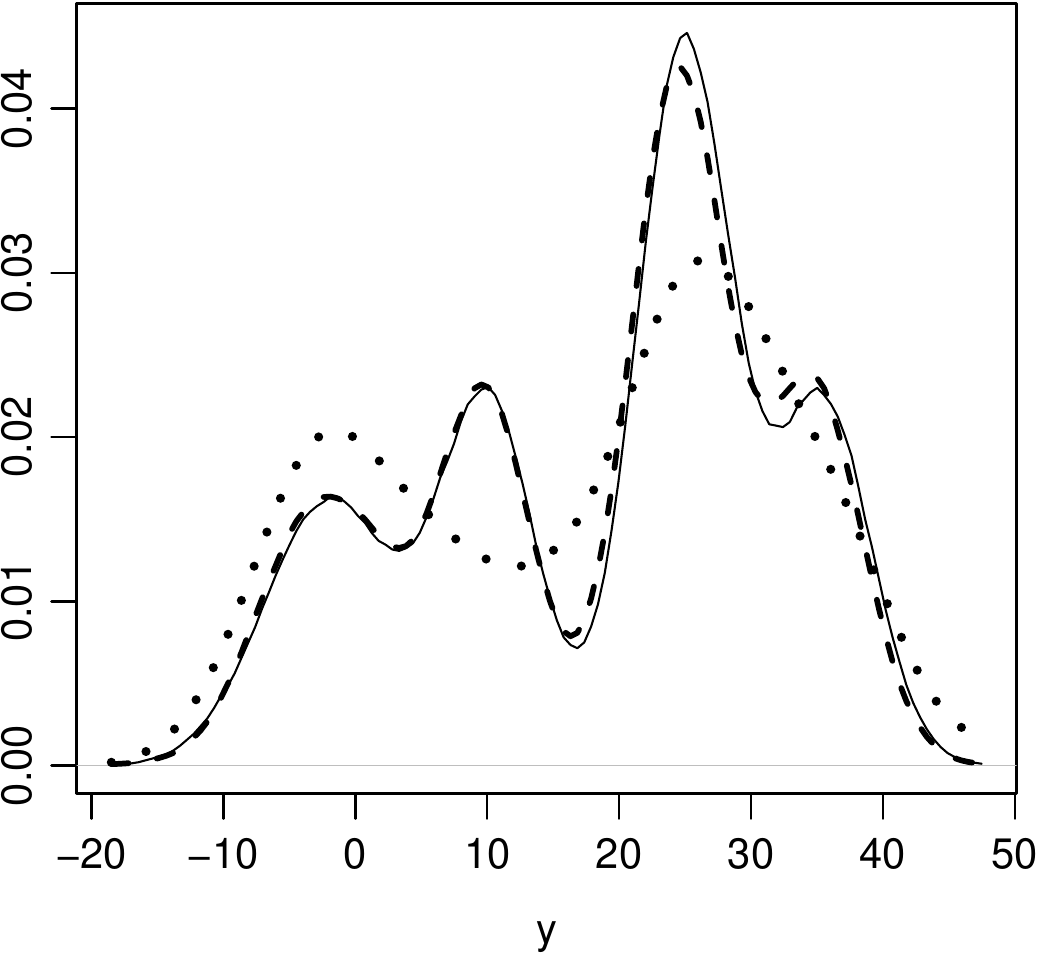}}
	\end{tabular}
	\caption{Typical realizations of estimates for the density of the response variable $Y$ in \eqref{eq:mult-reg-multimodal} given by the Rosenblatt–Parzen density estimator $\fhY^{(RP)}$ (dotted line), and the multiple regression-enhanced convolution estimator $\fhY^{(MR)}$ (dashed line). The black solid line is the true density $\fY$. \label{fig:typ-real-multimodal}}
\end{figure}

Next, we consider the case of a response variable with a multimodal distribution,
\begin{align} \label{eq:mult-reg-multimodal}
Y = \alpha_0 + \alpha_1 X + \eps[],
\end{align}
where $(\alpha_0,\alpha_1) = (4,1.5)$, $\eps[] \sim \mathcal{N}(0,4)$, and the density $\fX$ of $X$ is given by 
\begin{align*}
\fX(y) = \sum_{i=1}^4 w_i \Psi_i(y;\mu_i,\sigma_i),
\end{align*}
with $(w_1,w_2,w_3,w_4) = (0.2,0.2,0.4,0.2)$, $(\mu_1,\mu_2,\mu_3,\mu_4) = (-4,4,14,21)$, $(\sigma_1,\sigma_2,\sigma_3,\sigma_4) = (3,2,2,2)$. In Figure \ref{fig:typ-real-multimodal}, we plot $\fY$ along with typical realizations of density estimates given by $\fhY^{(RP)}$ and $\fhY^{(MR)}$, for a complete case dataset of size $N=100$, while the ratio of additional covariate observations supplied to the convolution estimator progressively increases, $\tau \in \{0,4,16,64\}$. The Rozenblatt-Parzen $\fhY^{(RP)}$ completely fails to resolve two of the modes of this distribution and severely underestimates the largest mode. At $\tau=0$, the convolution estimator $\fhY^{(MR)}$ manages to pick out three modes, although the magnitudes of these modes and the general shape of the density is not so accurate. However, as additional covariate observations are incorporated, the convolution estimator provides an increasingly accurate representation of the true density.

In Table \ref{tab:mult-reg-multimodal}, we report MISE results for both $\fhY^{(RP)}$ and $\fhY^{(MR)}$. Three complete cases samples sizes are considered, $N \in \{50,100,200\}$. Once again a very substantial substantial reduction in MISE is observed, with $\MISE[\fhY^{(MR)}]$ about $35$ times smaller than $\MISE[\fhY^{(RP)}]$ at $\tau = 128$ and $N=200$.

\renewcommand{\arraystretch}{1.4}
\setlength{\tabcolsep}{.4em}
\begin{table}[ht]
\footnotesize
\centering
\begin{tabular}[t]{|l|c|cccccccccc|}
\hline
      & $\MISE[\fhY^{(RP)}]$ & \multicolumn{10}{c|}{$\MISE[\fhY^{(MR)}]$} \\
      &                & $\tau=0$ & $\tau=2$ & $\tau=4$ & $\tau=8$ & $\tau=16$ & $\tau=32$ & $\tau=64$ & $\tau=128$ & $\tau=256$ & $\tau=512$ \\
\hline
$N=50$  &  $2.43\text{e-}03$               & $1.89\text{e-}03$ & $8.08\text{e-}04$ & $5.35\text{e-}04$ & $3.28\text{e-}04$ & $2.24\text{e-}04$ &   $1.59\text{e-}04$     & $1.22\text{e-}04$ & $1.03\text{e-}04$ & $8.76\text{e-}05$ & $8.83\text{e-}05$ \\
$N=100$ &  $1.54\text{e-}03$               & $1.05\text{e-}03$ & $4.19\text{e-}04$ & $2.64\text{e-}04$ & $1.71\text{e-}04$ & $1.09\text{e-}04$ &   $7.34\text{e-}05$     & $6.06\text{e-}05$ & $5.23\text{e-}05$ & $4.94\text{e-}05$ & $4.22\text{e-}05$ \\
$N=200$ &  $9.15\text{e-}04$               & $5.90\text{e-}04$ & $2.17\text{e-}04$ & $1.35\text{e-}04$ & $8.45\text{e-}05$ & $5.57\text{e-}05$ &   $3.92\text{e-}05$     & $2.98\text{e-}05$ & $2.70\text{e-}05$ & $2.29\text{e-}05$ & $2.37\text{e-}05$ \\
\hline
\end{tabular}
\caption{MISE results for the Rozenblatt-Parzen estimator and the multiple regression-enhanced convolution estimator for the density of the response variable in the regression model \eqref{eq:mult-reg-corr-cov}. \label{tab:mult-reg-multimodal}}
\end{table}

\subsection{Multiple regression with correlated covariates and non-Gaussian error}
\begin{figure}[ht!]
	\centering
	\captionsetup[subfloat]{labelformat=empty}
	\begin{tabular}{ccc}
	\subfloat[(i) $\tau=0$]{\includegraphics[scale=0.6]{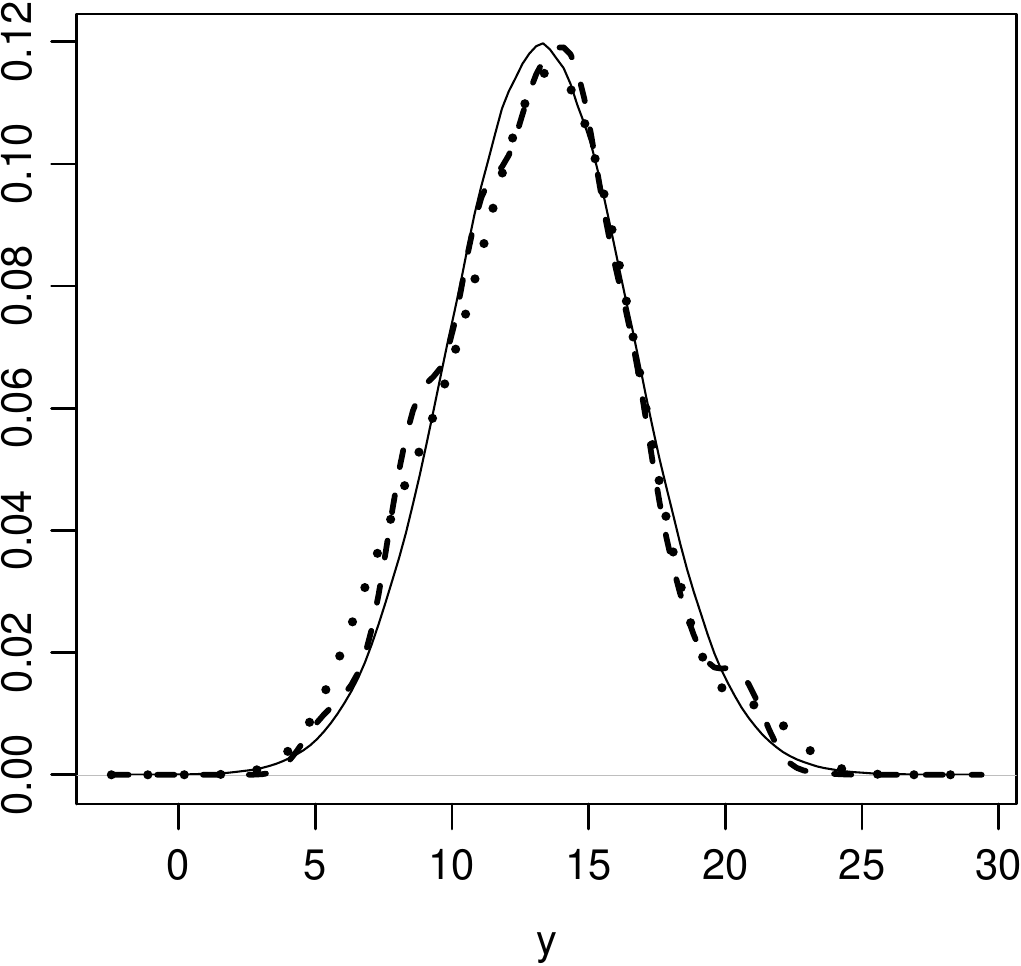}} \hspace{2em}
	\subfloat[(ii) $\tau=4$]{\includegraphics[scale=0.6]{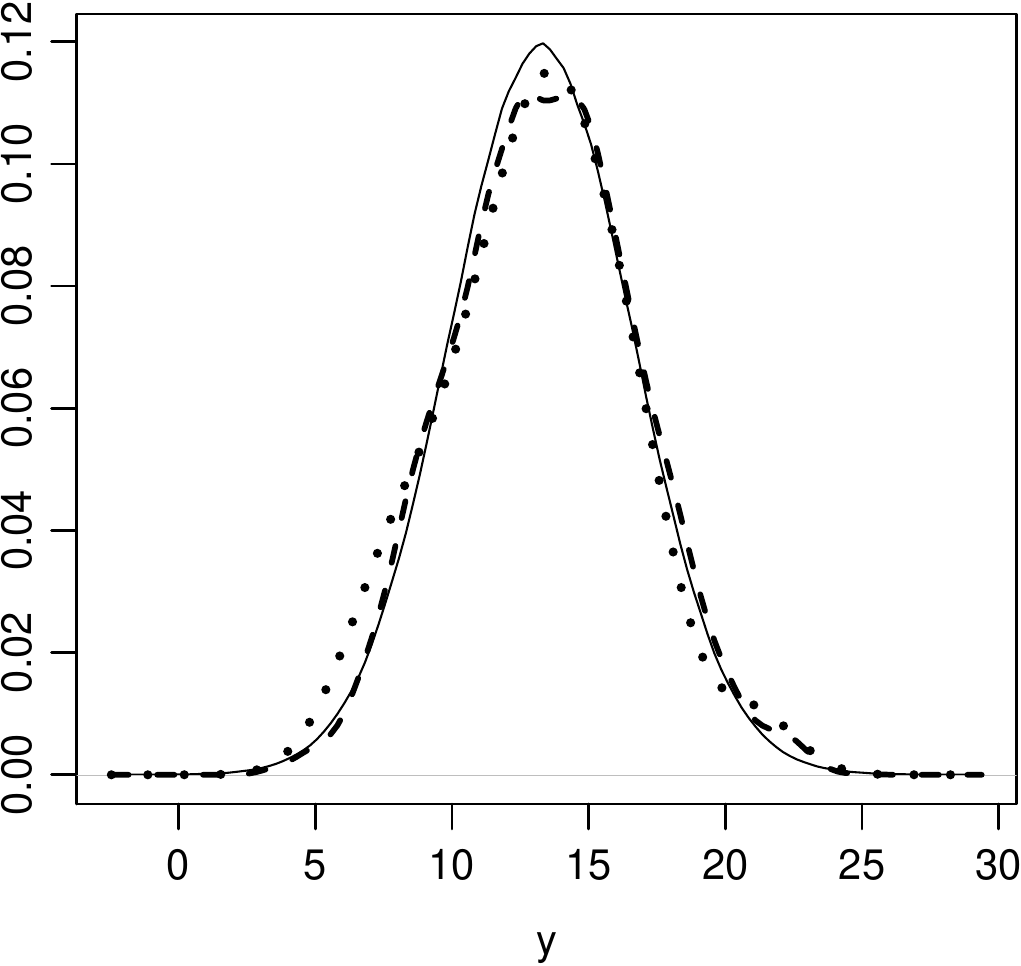}} \\
	\subfloat[(iii) $\tau=16$]{\includegraphics[scale=0.6]{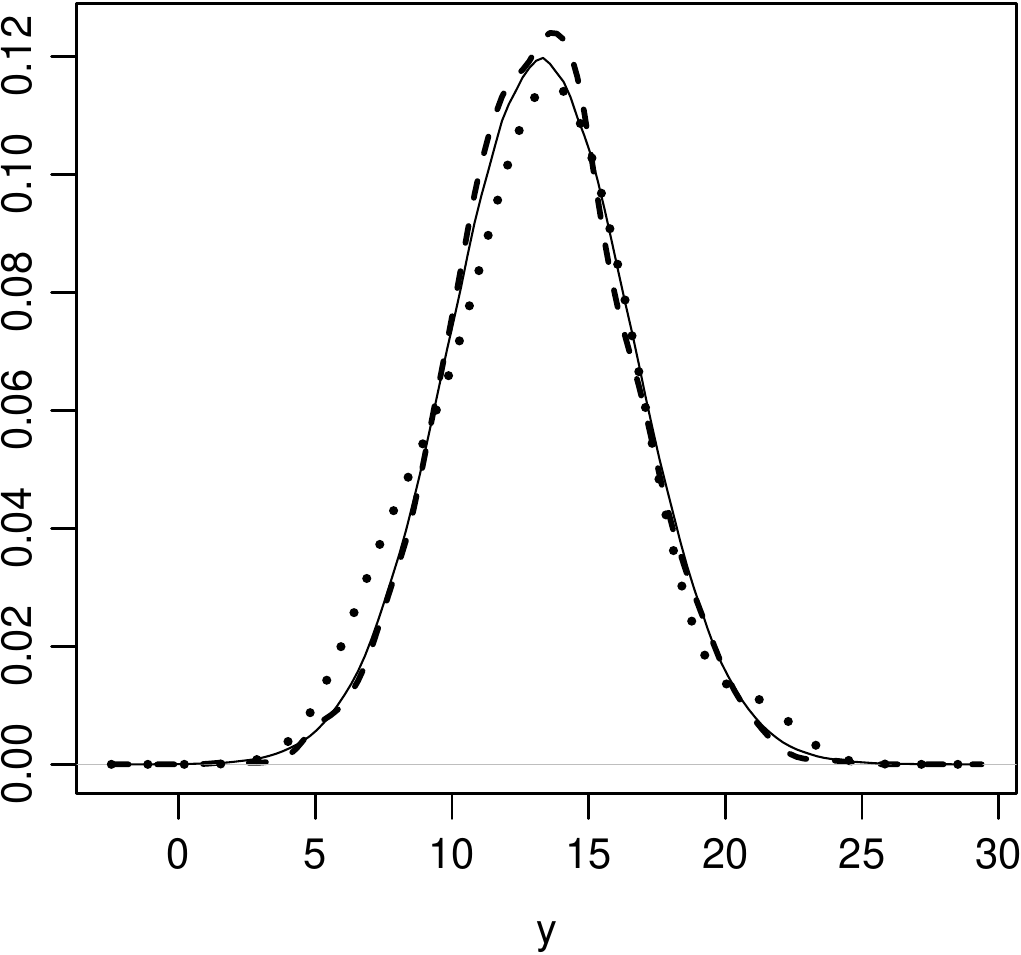}} \hspace{2em}
	\subfloat[(iv) $\tau=64$]{\includegraphics[scale=0.6]{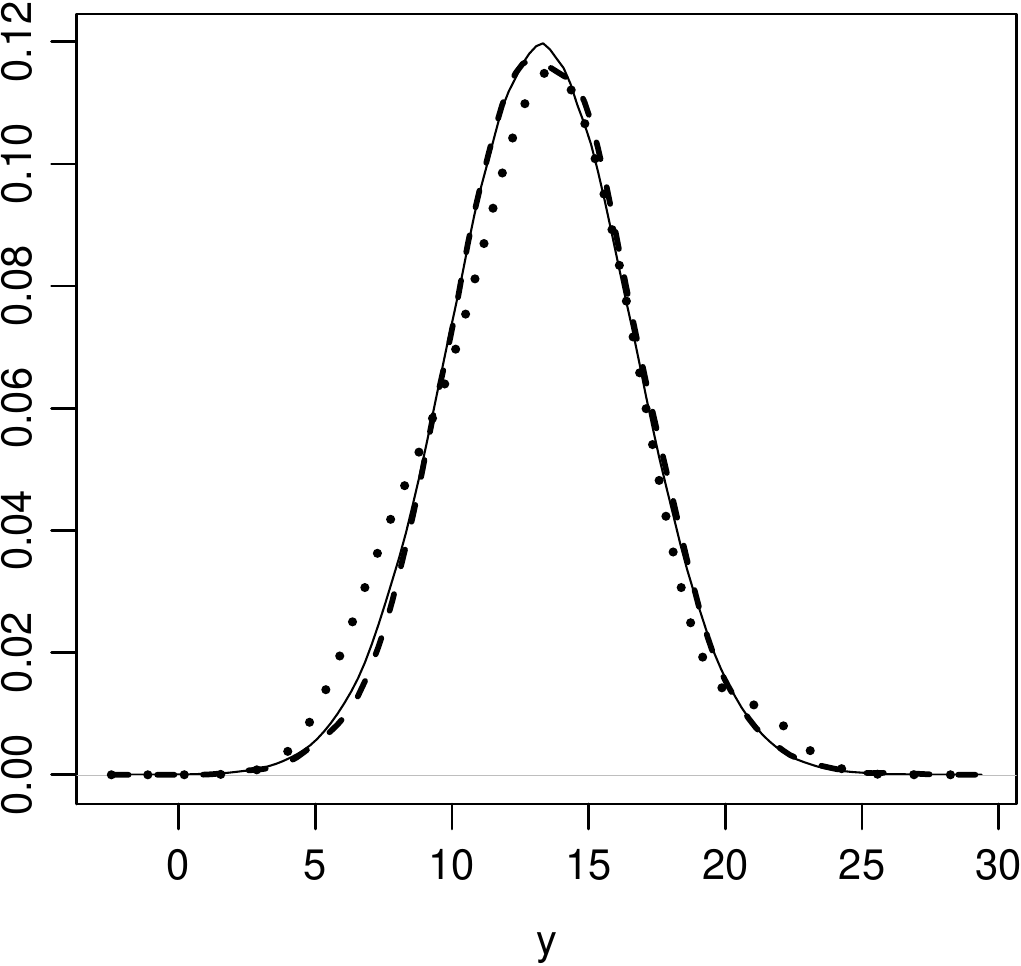}}
	\end{tabular}
	\caption{Typical realizations of estimates for the density of the response variable $Y$ in \eqref{eq:mult-reg-corr-cov} given by the Rosenblatt–Parzen density estimator $\fhY^{(RP)}$ (dotted line), and the multiple regression-enhanced convolution estimator $\fhY^{(MR)}$ (dashed line). The black solid line is the true density $\fY$. \label{fig:typ-real-corr-cov}}
\end{figure}

Next we consider a multiple regression model where the covariates are correlated and the error is non-Gaussian,
\begin{align} \label{eq:mult-reg-corr-cov}
Y = \alpha_0 + \alpha_1 X_1 + \alpha_2 X_2 + \alpha_3 X_3 + \eps[],
\end{align}
where $(\alpha_0,\alpha_1,\alpha_2,\alpha_3) = (1,1,2,0.5)$, the covariates are distributed as
\begin{align*}
X_1 \sim \beta(2,5), \quad \quad
X_2 \sim \mathcal{N}(6,4), \quad \quad X_3 \sim t_{6}, \quad \quad
\eps[] & \sim \text{Skew-Normal}(\xi,\omega,\alpha),
\end{align*}
and the correlation matrix for the covariates is defined as
\begin{align*}
\Corr(X)
&=
\begin{bmatrix}
1 & 0.2 & 0.5, \\
0.2 & 1 & 0.3, \\
0.5 & 0.3 & 1
\end{bmatrix}.
\end{align*}
The parameters for the skew normal error distribution are set to $(\xi,\omega,\alpha) = (-\omega  \sqrt{2 \alpha^2/((1+\alpha^2)\pi)},1,3)$, where $\xi$ has been chosen to ensure that the error has mean zero. In Figure \ref{fig:typ-real-corr-cov}, we plot $\fY$ along with typical realizations of density estimates given by $\fhY^{(RP)}$ and $\fhY^{(MR)}$.
%
MISE results for both estimators are presented in Table \ref{tab:corr-cov-non-neg-error}. This time $\MISE[\fhY^{(MR)}]$ is about $32$ times smaller than $\MISE[\fhY^{(RP)}]$ at $\tau = 128$ and $N=200$.

\renewcommand{\arraystretch}{1.4}
\setlength{\tabcolsep}{.4em}
\begin{table}[htb!]
\footnotesize
\centering
\begin{tabular}[t]{|l|c|cccccccccc|}
\hline
      & $\MISE[\fhY^{(RP)}]$ & \multicolumn{10}{c|}{$\MISE[\fhY^{(MR)}]$} \\
      &                & $\tau=0$ & $\tau=2$ & $\tau=4$ & $\tau=8$ & $\tau=16$ & $\tau=32$ & $\tau=64$ & $\tau=128$ & $\tau=256$ & $\tau=512$ \\
\hline
$N=50$  &  $3.66\text{e-}03$               & $4.94\text{e-}03$ & $1.87\text{e-}03$ & $1.24\text{e-}03$ & $7.73\text{e-}04$ & $4.51\text{e-}04$ &   $2.81\text{e-}04$     & $1.92\text{e-}04$ & $1.39\text{e-}04$ & $1.27\text{e-}04$ & $1.08\text{e-}04$ \\
$N=100$ &  $2.04\text{e-}03$               & $2.83\text{e-}03$ & $1.08\text{e-}03$ & $6.72\text{e-}04$ & $3.92\text{e-}04$ & $2.25\text{e-}04$ &   $1.41\text{e-}04$     & $9.56\text{e-}05$ & $7.04\text{e-}05$ & $5.55\text{e-}05$ & $5.45\text{e-}05$ \\
$N=200$ &  $1.15\text{e-}03$               & $1.51\text{e-}03$ & $5.81\text{e-}04$ & $3.58\text{e-}04$ & $2.02\text{e-}04$ & $1.17\text{e-}04$ &   $7.16\text{e-}05$     & $4.83\text{e-}05$ & $3.58\text{e-}05$ & $3.01\text{e-}05$ & $2.64\text{e-}05$ \\
\hline
\end{tabular}
\caption{MISE results for the Rozenblatt-Parzen estimator and the multiple regression-enhanced convolution estimator for the density of the response variable in the regression model \eqref{eq:mult-reg-corr-cov}. \label{tab:corr-cov-non-neg-error}}
\end{table}

\section{Concluding remarks} \label{sec:concluding-remarks}
In this work, we have proposed a convolution estimator for enhancing the accuracy of estimates of the density of a response variable in a sample of $N$ complete case observations, by using an additional sample of $M$ covariate observations. While previous works on convolution estimators have modelled the relationship between the response variable and the covariates using nonlinear regression models, in this paper a multiple regression model was employed. Unlike Nadaraya-Watson-based convolution estimators that suffer from the curse of dimensionality, we showed that the convergence of the multiple regression-enhanced convolution estimator is independent of the dimensionality of the covariates, which is due to the fact that the convergence of the underlying OLS estimator is also dimension independent.

The case of $M$ additional covariate observations is a generalization of the usual convolution estimator setting considered in the literature. The usual setting involves estimating the density of a response variable using a sample of $N$ complete case observations of a response variable and an associated set of covariates. By setting $M = 0$, we recover this case. Indeed, we showed that $\MSE[\fhY(y)] = O(h^4) + O((hNL)^{-1}) + O(L^{-1}) + O(N^{-1})$ reduces to $\MSE[\fhY(y)] = O(h^4) + O((hN^2)^{-1}) + O(N^{-1})$ when $M=0$, which recovers previous results in the literature \cite[Eq. (23)]{stove2012convolution}and \cite[Eq. (3.4)]{escanciano2012n}.

By deriving the asymptotic MSE and the asymptotically optimal bandwidth, we resolved the question on the convergence of convolution estimators with respect to the size of the additional sample of $M$ covariate observations. That is, we proved that for a large fixed $N$, at the asymptotically optimal bandwidth, the MSE converges as $O(M^{-4/5})$ towards an $O(N^{-1})$ constant. We also showed that for a fixed $M$, the MSE converges as $O(N^{-1})$. This means that supplying the convolution estimator with additional covariate observations is not quite as effective as supplying it with more complete cases observations. Crucially, however, in many practical applications it can be difficult if not downright impossible to obtain more observations of a response variable, while at the same time it can be very straightforward to obtain more observations of the covariates. 

Numerical simulations confirmed the existence of the saturation phenomena predicted by the theory, whereby the MISE converges as $O(M^{-4/5})$ towards an $O(N^{-1})$ constant, as opposed to zero, as the size of the additional sample increases, whereas it converges as $O(N^{-1})$ towards zero as the size of the complete case sample increases. Moreover numerical simulations demonstrated that even if the MISE of the multiple regression-enhanced convolution estimator is greater than that of the Rosenblatt-Parzen density estimator on the complete case dataset, by supplying the convolution estimator with additional covariate observations, its MISE can be made about $20$ to $35$ times smaller than the MISE of the Rosenblatt-Parzen before the accuracy improvement saturates.

The evaluation of the multiple regression-enhanced convolution estimator is an order of magnitude more computationally expensive than the evaluation of the Rosenblatt–Parzen density estimator. To reduce computational costs, we developed a FGT-based acceleration algorithm that draws on the high performance C++ library FIGTree. Simulations showed that this algorithm dramatically outperformed a variety of alternative evaluation approaches. In particular, it was demonstrated to be up to $40$ times faster than FFT-based acceleration, with the reduction in computational times becoming even more pronounced as $M$ increases.

In terms of future research directions, some interesting topics include heteroscedasticity, segmented regression-enhanced convolution estimators, and transfer learning for cases when the $N$ complete case observations and $M$ additional covariate observations have different distributions. Convolution estimators in the presence of heteroscedasticity have been considered in works such as \cite{stove2012convolution} and \cite{li2016n}. As we have demonstrated in this paper, however, to obtain fast convergence rates of $O(N^{-1})$ and $O(M^{-4/5})$ that are independent of the covariate dimensionality, the underlying regression function estimator has to be immune to the curse of dimensionality. Thus, we anticipate the method chosen to model heteroscedastic errors also needs to be immune to the curse of dimensionality if one wants to main these fast convergence rates.

In cases where the data is not well fit by a multiple regression model, it might still be possible to fit it with a piecewise multiple regression model, as opposed to employing a fully nonlinear regression model. For example, in situations where the data is well fit by a multiple regression model in several segments, each with a different, albeit constant, error variance, a multidimensional segmented regression model could be employed \cite{diakonikolas2020efficient,liu1997segmented}. Since segmented regression involves partitioning the data into several segments, each with its own dedicated OLS estimator, we conjecture that a segmented regression-enhanced convolution estimator could retain the fast dimension independent convergence rates of the multiple regression-enhanced convolution estimator. Of course, segmented regression itself can be subject to the curse of dimensionality if too many segments are used, but this can be avoided by placing an upper bound on the number of partitions \cite{liu1997segmented}.

Transfer learning is another very interesting avenue for future research with regards to convolution estimators. To the best of our knowledge this topic is completely unexplored. Transfer learning has been demonstrated to be very effective at utilizing labelled information from a source domain to enhance the performance of a model in a separate target domain with little or no labelled data \cite{day2017survey,zhuang2020comprehensive,pan2009survey}. Our work in this paper was concerned with the case where the $N$ complete case observations and the $M$ additional covariate observations came from the same distribution. There are many practical applications in which the complete case observations and additional covariate observations could come from different, yet closely, related distributions. We expect that incorporating transfer learning capabilities into convolution estimators could provide a significant improvement in accuracy in situations such as these.

\section*{}

\bibliography{mybibfile}

\begin{thebibliography}{10}
\expandafter\ifx\csname url\endcsname\relax
  \def\url#1{\texttt{#1}}\fi
\expandafter\ifx\csname urlprefix\endcsname\relax\def\urlprefix{URL }\fi
\expandafter\ifx\csname href\endcsname\relax
  \def\href#1#2{#2} \def\path#1{#1}\fi

\bibitem{rosenblatt1956}
M.~Rosenblatt, \href{https://doi.org/10.1214/aoms/1177728190}{Remarks on some
  nonparametric estimates of a density function}, Ann. Math. Statist. 27~(3)
  (1956) 832--837.
\newblock \href {http://dx.doi.org/10.1214/aoms/1177728190}
  {\path{doi:10.1214/aoms/1177728190}}.
\newline\urlprefix\url{https://doi.org/10.1214/aoms/1177728190}

\bibitem{parzen1962estimation}
E.~Parzen, On estimation of a probability density function and mode, The annals
  of mathematical statistics 33~(3) (1962) 1065--1076.

\bibitem{escanciano2012n}
J.~C. Escanciano, D.~T. Jacho-Ch{\'a}vez, n-uniformly consistent density
  estimation in nonparametric regression models, Journal of Econometrics
  167~(2) (2012) 305--316.

\bibitem{muller2012estimating}
U.~U. M{\"u}ller, Estimating the density of a possibly missing response
  variable in nonlinear regression, Journal of Statistical Planning and
  Inference 142~(5) (2012) 1198--1214.

\bibitem{stove2012convolution}
B.~St{\o}ve, D.~Tj{\o}stheim, A convolution estimator for the density of
  nonlinear regression observations, Scandinavian Journal of Statistics 39~(2)
  (2012) 282--304.

\bibitem{li2016n}
S.~Li, Y.~Tu, n-consistent density estimation in semiparametric regression
  models, Computational Statistics \& Data Analysis 104 (2016) 91--109.

\bibitem{schick2004root}
A.~Schick*, W.~Wefelmeyer, Root n consistent density estimators for sums of
  independent random variables, Journal of Nonparametric Statistics 16~(6)
  (2004) 925--935.

\bibitem{schick2007root}
A.~Schick, W.~Wefelmeyer, Root-n consistent density estimators of convolutions
  in weighted l1-norms, Journal of Statistical Planning and Inference 137~(6)
  (2007) 1765--1774.

\bibitem{saavedra1999rate}
{\'A}.~Saavedra, R.~Cao, Rate of convergence of a convolution-type estimator of
  the marginal density of a ma (1) process, Stochastic processes and their
  applications 80~(2) (1999) 129--155.

\bibitem{saavedra2000estimation}
A.~Saavedra, R.~Cao, On the estimation of the marginal density of a moving
  average process, Canadian Journal of Statistics 28~(4) (2000) 799--815.

\bibitem{gyorfi2006distribution}
L.~Gy{\"o}rfi, M.~Kohler, A.~Krzyzak, H.~Walk, A distribution-free theory of
  nonparametric regression, Springer Science \& Business Media, 2006.

\bibitem{morariu2008automatic}
V.~Morariu, B.~Srinivasan, V.~C. Raykar, R.~Duraiswami, L.~S. Davis, Automatic
  online tuning for fast gaussian summation, Advances in neural information
  processing systems 21 (2008) 1113--1120.

\bibitem{greengard1991fast}
L.~Greengard, J.~Strain, The fast gauss transform, SIAM Journal on Scientific
  and Statistical Computing 12~(1) (1991) 79--94.

\bibitem{arya1993approximate}
S.~Arya, D.~M. Mount, Approximate nearest neighbor queries in fixed
  dimensions., in: SODA, Vol.~93, Citeseer, 1993, pp. 271--280.

\bibitem{greene2003econometric}
W.~H. Greene, Econometric analysis, Pearson Education India, 2003.

\bibitem{silverman1986density}
B.~W. Silverman, Density estimation for statistics and data analysis, Vol.~26,
  CRC press, 1986.

\bibitem{rudemo1982empirical}
M.~Rudemo, Empirical choice of histograms and kernel density estimators,
  Scandinavian Journal of Statistics (1982) 65--78.

\bibitem{sheather1991reliable}
S.~J. Sheather, M.~C. Jones, A reliable data-based bandwidth selection method
  for kernel density estimation, Journal of the Royal Statistical Society:
  Series B (Methodological) 53~(3) (1991) 683--690.

\bibitem{deng2011density}
H.~Deng, H.~Wickham, Density estimation in r, Electronic publication.

\bibitem{diakonikolas2020efficient}
I.~Diakonikolas, J.~Li, A.~Voloshinov, Efficient algorithms for
  multidimensional segmented regression, arXiv preprint arXiv:2003.11086.

\bibitem{liu1997segmented}
J.~Liu, S.~Wu, J.~V. Zidek, On segmented multivariate regression, Statistica
  Sinica (1997) 497--525.

\bibitem{day2017survey}
O.~Day, T.~M. Khoshgoftaar, A survey on heterogeneous transfer learning,
  Journal of Big Data 4~(1) (2017) 1--42.

\bibitem{zhuang2020comprehensive}
F.~Zhuang, Z.~Qi, K.~Duan, D.~Xi, Y.~Zhu, H.~Zhu, H.~Xiong, Q.~He, A
  comprehensive survey on transfer learning, Proceedings of the IEEE 109~(1)
  (2020) 43--76.

\bibitem{pan2009survey}
S.~J. Pan, Q.~Yang, A survey on transfer learning, IEEE Transactions on
  knowledge and data engineering 22~(10) (2009) 1345--1359.

\bibitem{gut2013probability}
A.~Gut, Probability: a graduate course, Vol.~75, Springer Science \& Business
  Media, 2013.

\end{thebibliography}

\appendix
\begin{appendices}
\section{} \label{appendix:asy-expectations}

Our approach in this section is quite similar to that of St{\o}ve and Tj{\o}stheim in \cite[Supp. Material]{stove2012convolution}. In fact, we use the same approximation these authors introduce in \cite[Supp. Material, Proof of Theorem 2]{stove2012convolution}; suppose $n$ rows, given by the index set $\Ib = \{i_j\}_{j=1}^n$, where $i_j \in \{1,2,\dots,N\}$, are removed from the multiple regression model \eqref{eq:samnple-mult-reg}, with $n \in \{1,2,3,4\}$. Denote by $\alphabh[I]$ the OLS estimator associated with this reduced multiple regression model.

We use the reduced OLS estimator to approximate $K_{ij}(\alphabh[])$ by $K_{ij}(\alphabh[I])$ in the expectations in this section, where we recall that $K_{ij}$ is given by \eqref{eq:K-ij}. Intuitively speaking, this change can be ignored asymptotically because removing a very small finite number of rows from the multiple regression model \eqref{eq:samnple-mult-reg} has an asymptotically negligible effect on the convergence as the sample size $N \to \infty$.

Denote by
\begin{align*}
f_{\epsb[i_1 i_2 \dots i_n-]}(\eb[i_1 i_2 \dots i_n-])
= \prod_{\substack{i=1 \\ i \notin \{i_1, i_2, \dots, i_n\} }}^N f_{\epsb[i]}(\eb[i]), \quad \quad  \quad \quad
d\eb[i_1, i_2, \dots, i_n-]= \prod_{\substack{i=1 \\ i \notin \{i_1, i_2, \dots, i_n\} }}^N d\eb[i].
\end{align*}
This notation allows us to write expressions such as
\begin{align*}
\int g(e_1,e_2,\dots,e_N) \prod_{i=1}^N f_{\eps[i]}(e_i) \prod_{i=1}^N d e_i,
\end{align*}
in a form that is more convenient for the analysis in this section, namely,
\begin{align*}
\int g(e_1,e_2,\dots,e_N) f_{\epsb[12-]}(\eb[12-])f_{\eps[1]}(e_1)f_{\eps[2]}(e_2) d\eb[12-] de_1 de_2,
\end{align*}
where $g$ is some arbitrary function.  The lemmas in this section hold under assumptions (A), (B), (C), (D), (E), and (F) given in Section \ref{subsec:Assumptions}.

%
%
%
%
%
%
%
%
\begin{lemma}[] \label{lem:E-K-11-leading-order}
As $h \to 0$, it holds that
\begin{align*}
E[K_h(y - \Xb[1]^T\alphab[] - \eps[1])] \sim h \fY(y) + h^3 \frac{\mu_K}{2} \fY''(y).
\end{align*}
\begin{proof}
Performing the change of variables $-r = (y - \xb[1]^T \alphab[] - \eps[1])/h$, using the fact that $K$ is a symmetric function, Taylor expanding with respect to $h$, and making use of \eqref{eq:int-K}, \eqref{eq:int-K-defs}, and \eqref{eq:fYk-conv-form}, we find that
\begin{align*}
& E[K_h(y - \Xb[1]^T \alphab[] - \eps[1])] \\
& = \int \cdots \int K_h(y - \xb[1]^T \alphab[] - e_1) f_{\eps[1]}(e_1) f_{\Xb[1]}(\xb[1]) \ d\xb[1] de_1 \\
& = h \int \cdots \int K(r) f_{\eps[1]}(y - \xb[1]^T \alphab[] + hr) f_{\Xb[1]}(\xb[1]) \ d\xb[1] dr \\
& \sim h \int \cdots \int K(r) (f_{\eps[1]}(y - \xb[1]^T \alphab[]) + hr f_{\eps[1]}'(y - \xb[1]^T \alphab[]) + \frac{(hr)^2}{2}f_{\eps[1]}''(y - \xb[1]^T \alphab[])) f_{\Xb[1]}(\xb[1]) \ d\xb[1] dr \\
& = h \int (f_{\eps[1]}(y - \xb[1]^T \alphab[]) + h^2 \frac{\mu_K}{2}f_{\eps[1]}''(y - \xb[1]^T \alphab[])) f_{\Xb[1]}(\xb[1]) \ d\xb[1] \\
& = h E[f_{\eps[]}(y - \Xb[1]^T \alphab[])] + h^3 \frac{\mu_K}{2} E[f_{\eps[]}''(y - \Xb[1]^T \alphab[])] \\
& = h \fY(y) + h^3 \frac{\mu_K}{2} \fY''(y).
\end{align*}
\end{proof}
\end{lemma}

%
%
%
%
%
%
%
%
\begin{lemma} \label{lem:E-feps-2}
As $N \to \infty$, it holds that
\begin{align*}
E[f_{\eps[]}(\yt(\alphabh[I],\Xb[1],\Xb[2]))]
& \sim \fY(y)
+ N^{-1} \frac{\sigmaeps^2}{2} \sum_{p_1,p_2=0}^J E[(\Phii[N]^{-1} \Xb[1])_{p_1}(\Phii[N]^{-1} \Xb[1])_{p_2}] \\
& \times E[f_{\eps[]}''(y - \Xb[1]^T \alphab[]) (\Xb[2] - \Xb[1])_{p_1}(\Xb[2] - \Xb[1])_{p_2}].
\end{align*}
\begin{proof}
Taylor expanding $f_{\eps[]}(\yt(\alphabh[I],\Xb[1],\Xb[2]))$ with respect to $\alphabh[I]$ about $\alphab[]$, taking the expectation, and using \eqref{eq:yt-simple}, we get
\begin{align*}
E[f_{\eps[]}(\yt(\alphabh[I],\Xb[1],\Xb[2]))]
&= E[f_{\eps[]}(y - \Xb[1]^T \alphab[])]
+ \sum_{p_1=0}^J E[(\alphabh[I]-\alphab[])_{p_1}] E\bigg[\frac{\partial f_{\eps[]}(\yt(\alphabh[I],\Xb[1],\Xb[2]))}{\partial(\alphabh[I])_{p_1}} \bigg|_{\alphabh[I]=\alphab[]}\bigg] \\
& + \frac{1}{2} \sum_{p_1,p_2=0}^J E[(\alphabh[I]-\alphab[])_{p_1}(\alphabh[I]-\alphab[])_{p_2}] E\bigg[\frac{\partial^2 f_{\eps[]}(\yt(\alphabh[I],\Xb[1],\Xb[2]))}{\partial(\alphabh[I])_{p_1}\partial(\alphabh[I])_{p_2}} \bigg|_{\alphabh[I]=\alphab[]}\bigg] \\
& + \frac{1}{6} \sum_{p_1,p_2,p_3=0}^J E[(\alphabh[I]-\alphab[])_{p_1}(\alphabh[I]-\alphab[])_{p_2}(\alphabh[I]-\alphab[])_{p_3} \frac{\partial^3 f_{\eps[]}(\yt(\alphabh[I],\Xb[1],\Xb[2]))}{\partial(\alphabh[I])_{p_1}\partial(\alphabh[I])_{p_2}\partial(\alphabh[I])_{p_3}} \bigg|_{\alphabh[I]=\zetab[]}\bigg],
\end{align*}
where $\zetab[] = \alphab[] + c(\alphabh[I] - \alphab[])$ for $c \in (0,1)$. The first order term vanishes since $E[(\alphabh[I]-\alphab[])_{p_1}]=0$. Note that once the derivatives in the first and second order terms are evaluated at $\alphab[I] = \alphab[]$ they become independent of $\alphab[I]$ since they only depend on $\Xb[1]$ and $\Xb[2]$, and these covariate observations are not present in $\alphab[I]$. On the other hand, in the remainder term, $\zetab[]$ does depend on $\alphab[I]$. However, at leading order as $N \to \infty$, $\alphabh[I] \sim \alphab[]$, which means that $\zetab[] \sim \alphab[]$. Thus,
\begin{align*}
E[f_{\eps[]}(\yt(\alphabh[I],\Xb[1],\Xb[2]))]
&\sim \fY(y)
+ \frac{1}{2} \sum_{p_1,p_2=0}^J E[(\alphabh[I]-\alphab[])_{p_1}(\alphabh[I]-\alphab[])_{p_2}] E\bigg[\frac{\partial^2 f_{\eps[]}(\yt(\alphabh[I],\Xb[1],\Xb[2]))}{\partial(\alphabh[I])_{p_1}\partial(\alphabh[I])_{p_2}} \bigg|_{\alphabh[I]=\alphab[]}\bigg] \\
& + \frac{1}{6} \sum_{p_1,p_2,p_3=0}^J E[(\alphabh[I]-\alphab[])_{p_1}(\alphabh[I]-\alphab[])_{p_2}(\alphabh[I]-\alphab[])_{p_3}] E\bigg[\frac{\partial^3 f_{\eps[]}(\yt(\alphabh[I],\Xb[1],\Xb[2]))}{\partial(\alphabh[I])_{p_1}\partial(\alphabh[I])_{p_2}\partial(\alphabh[I])_{p_3}} \bigg|_{\alphabh[I]=\alphab[]}\bigg],
\end{align*}
where we used \eqref{eq:fY-conv-form} for the leading-order term. Next, we can approximate $\alphabh[I]$ by $\alphabh[]$ as $N\to \infty$, and then use \eqref{eq:OLS-estimator-Phi} and the independence of the error observations to get
\begin{align*}
E[(\alphabh[I]-\alphab[])_{p_1}(\alphabh[I]-\alphab[])_{p_2}]
& \sim E[(\alphabh[]-\alphab[])_{p_1}(\alphabh[]-\alphab[])_{p_2}] \\
&= E\bigg[(N^{-1} \Phii[N]^{-1} \sum_{_1i=1}^N \Xb[i_1]\eps[i_1])_{p_1}(N^{-1} \Phii[N]^{-1} \sum_{i_2=1}^N \Xb[i_2]\eps[i_2])_{p_2}\bigg] \\
&\sim N^{-2} \sum_{i_1,i_2=1}^N E[\eps[i_1]\eps[i_2]]E[(\Phii[N]^{-1}  \Xb[i_1])_{p_1}(\Phii[N]^{-1} \Xb[i_2])_{p_2}] \\
&= N^{-2} \sum_{i=1}^N E[\eps[i]\eps[i]] E[(\Phii[N]^{-1}  \Xb[i])_{p_1}(\Phii[N]^{-1} \Xb[i])_{p_2}] \\
&= N^{-1} \sigmaeps^2 E[(\Phii[N]^{-1}  \Xb[1])_{p_1}(\Phii[N]^{-1} \Xb[1])_{p_2}].
\end{align*}
Next, evaluating the derivative, and then retaining the leading-order term gives
\begin{align*}
E\bigg[\frac{\partial^2 f_{\eps[]}(\yt(\alphabh[I],\Xb[1],\Xb[2]))}{\partial(\alphabh[I])_{p_1}\partial(\alphabh[I])_{p_2}}\bigg|_{\alphabh[]=\alphab[]}\bigg]
& \sim E[f_{\eps[]}''(y - \Xb[1]^T \alphab[I]) (\Xb[2] - \Xb[1])_{p_1}(\Xb[2] - \Xb[1])_{p_2}] \\
& \sim E[f_{\eps[]}''(y - \Xb[1]^T \alphab[]) (\Xb[2] - \Xb[1])_{p_1}(\Xb[2] - \Xb[1])_{p_2}].
\end{align*}
Finally,
\begin{align*}
E[(\alphabh[I]-\alphab[])_{p_1}(\alphabh[I]-\alphab[])_{p_2}(\alphabh[I]-\alphab[])_{p_3}]
&\sim E[(\alphabh[]-\alphab[])_{p_1}(\alphabh[]-\alphab[])_{p_2}(\alphabh[]-\alphab[])_{p_3}] \\
&= N^{-3} \sum_{i_1,i_2,i_3=1}^N E[\eps[i_1]\eps[i_2]\eps[i_3]] E[(\Phii[N]^{-1}\Xb[i_1])_{p_1}(\Phii[N]^{-1} \Xb[i_2])_{p_2}(\Phii[N]^{-1} \Xb[i_3])_{p_3}] \\
&\sim N^{-3} \sum_{i=1}^N E[\eps[i]\eps[i]\eps[i]] E[(\Phii[N]^{-1}  \Xb[i])_{p_1}(\Phii[N]^{-1} \Xb[i])_{p_2}(\Phii[N]^{-1} \Xb[i])_{p_3}] \\
&= N^{-2} E[\eps[]^3] E[(\Phii[N]^{-1}  \Xb[1])_{p_1}(\Phii[N]^{-1} \Xb[1])_{p_2}(\Phii[N]^{-1} \Xb[1])_{p_3}] \\
&= O(N^{-2}),
\end{align*}
where we used Lemma \ref{lem:E-prod-PhiiN-Xi} for the last equality, which shows that the remainder is controlled.
\end{proof}
\end{lemma}

%
%
%
%
%
%
%
%
\begin{lemma} \label{lem:E-K-12}
As $N \to \infty$ and $h\to 0$, it holds that
\begin{align*}
E[K_{12}(\alphabh[])]
& \sim h \fY(y) + h^3 \frac{\mu_K}{2} \fY''(y) \\
& + h N^{-1} \frac{\sigmaeps^2}{2} \sum_{p_1,p_2=0}^J E[(\Phii[N]^{-1} \Xb[1])_{p_1}(\Phii[N]^{-1} \Xb[1])_{p_2}] E[f_{\eps[]}''(y - \Xb[1]^T \alphab[]) (\Xb[2] - \Xb[1])_{p_1}(\Xb[2] - \Xb[1])_{p_2}].
\end{align*}
\begin{proof}
Setting $\Ib = \{1,2\}$, we have that at leading-order, 
\begin{align*}
E[K_{12}(\alphabh[])]
\sim E[K_{12}(\alphabh[I])]
& = \int \cdots \int K_h(\yt(\alphabh[I],\xb[1],\xb[2]) - e_2) f_{\eps[2]}(e_2) f_{\epsb[12-]}(\eb[12-]) f_{\Xb[]}(\xb[]) \ d\xb[] d\eb[12-] de_2. 
\end{align*}
Performing the change of variables $-r = \yt(\alphabh[I],\xb[1],\xb[2]) - e_2)/h$, and using the fact that $K$ is a symmetric function, we have
\begin{align*}
E[K_{12}(\alphabh[])]
& \sim h \int \cdots \int K(r) f_{\eps[2]}(\yt(\alphabh[I],\xb[1],\xb[2]) + hr) f_{\epsb[12-]}(\eb[12-]) f_{\Xb[]}(\xb[]) \ d\xb[] d\eb[12-] dr.
\end{align*}
Next, Taylor expanding with respect to $h$ about $0$, and using \eqref{eq:int-K} and \eqref{eq:int-K-defs} gives
\begin{align*}
E[K_{12}(\alphabh[])]
& \sim h \int \cdots \int (f_{\eps[2]}(\yt(\alphabh[I],\xb[1],\xb[2])) + h^2 \frac{\mu_K}{2} f_{\eps[2]}''(\yt(\alphabh[I],\xb[1],\xb[2])) f_{\epsb[2-]}(\eb[2-]) f_{\Xb[]}(\xb[]) \ d\xb[] d\eb[12-] \\
& \sim h E[f_{\eps[]}(\yt(\alphabh[I],\Xb[1],\Xb[2]))] + h^3 \frac{\mu_K}{2} E[f_{\eps[]}''(\yt(\alphabh[I],\Xb[1],\Xb[2]))].
\end{align*}
The expression for the first expectation on the right hand side is provided in Lemma \ref{lem:E-feps-2}. On the other hand, taking the leading-order approximation $\alphab[]$ of $\alphabh[I]$, for the other expectation, we have that
\begin{align*}
E[f_{\eps[]}''(\yt(\alphabh[I],\Xb[1],\Xb[2]))]
= E[f_{\eps[]}''(y - \Xb[1]^T \alphab[] + (\Xb[2] - \Xb[1])^T(\alphabh[I] - \alphab[]))]
\sim E[f_{\eps[]}''(y - \Xb[1] \alphab[])] = \fY''(y),
\end{align*}
where we used \eqref{eq:fYk-conv-form} for the last equality.
\end{proof}
\end{lemma}

%
%
%
%
%
%
%
%
\begin{corollary} \label{cor:E-K-12-sqr}
As $N \to \infty$ and $h\to 0$, it holds that
\begin{align*}
E[K_{12}(\alphabh[])]^2
& \sim h^2 \fY^2(y) + h^4 \mu_K \fY(y) \fY''(y) \\
& + h^2 N^{-1} \sigmaeps^2  \fY(y) \sum_{p_1,p_2=0}^J E[(\Phii[N]^{-1} \Xb[1])_{p_1}(\Phii[N]^{-1} \Xb[1])_{p_2}] E[f_{\eps[]}''(y - \Xb[1]^T \alphab[]) (\Xb[2] - \Xb[1])_{p_1}(\Xb[2] - \Xb[1])_{p_2}].
\end{align*}
\end{corollary}

%
%
%
%
%
%
%
%
\begin{lemma} \label{lem:E-K-12-sqrt-internal}
As $N \to \infty$ and $h\to 0$, it holds that
\begin{align*}
E[K_{12}^2(\alphabh[])]
& \sim h \sigma_{K} \fY(y) + h^3 \frac{\sigma_{K,2}}{2} \fY''(y) \\
& + hN^{-1} \frac{\sigmaeps^2 \sigma_K}{2} \sum_{p_1,p_2=0}^J E[(\Phii[N]^{-1} \Xb[1])_{p_1}(\Phii[N]^{-1} \Xb[1])_{p_2}] E[f_{\eps[]}''(y - \Xb[1]^T \alphab[]) (\Xb[2] - \Xb[1])_{p_1}(\Xb[2] - \Xb[1])_{p_2}].
\end{align*}
\begin{proof}
Setting $\Ib = \{1,2\}$, we have that at leading-order
\begin{align*}
E[K_{12}^2(\alphabh[])]
\sim E[K_{12}^2(\alphabh[I])]
& = \int \cdots \int K_h^2(\yt(\alphabh[I],\xb[1],\xb[2]) - e_2) f_{\eps[2]}(e_2) f_{\epsb[2-]}(\eb[2-]) f_{\Xb[]}(\xb[]) \ d\xb[] d\eb[12-] de_2. 
\end{align*}
Performing the change of variables $-r = \yt(\alphabh[I],\xb[1],\xb[2]) - e_2)/h$, and using the fact that $K$ is a symmetric function, we have
\begin{align*}
E[K_{12}^2(\alphabh[])]
& \sim h \int \cdots \int K^2(r) f_{\eps[2]}(\yt(\alphabh[I],\xb[1],\xb[2]) + hr) f_{\epsb[2-]}(\eb[2-]) f_{\Xb[]}(\xb[]) \ d\xb[] d\eb[12-] dr.
\end{align*}
Next, Taylor expanding with respect to $h$ about $0$, and using \eqref{eq:int-K} and \eqref{eq:int-K-defs} gives
\begin{align*}
E[K_{12}^2(\alphabh[])]
& \sim h \sigma_{K} \int \cdots \int (f_{\eps[2]}(\yt(\alphabh[I],\xb[1],\xb[2])) + h^2 \mu_K \frac{1}{2} f_{\eps[2]}''(\yt(\alphabh[I],\xb[1],\xb[2])) f_{\epsb[2-]}(\eb[2-]) f_{\Xb[]}(\xb[]) \ d\xb[] d\eb[12-] \\
& \sim h \sigma_{K} E[f_{\eps[]}(\yt(\alphabh[I],\Xb[1],\Xb[2]))] + h^3 \frac{\sigma_{K,2}}{2} E[f_{\eps[]}''(\yt(\alphabh[I],\Xb[1],\Xb[2]))].
\end{align*}
Leading-order expressions for these expectations were already derived in Lemma \ref{lem:E-K-12}.
\end{proof}
\end{lemma}

%
%
%
%
%
%
%
%
\begin{lemma} \label{lem:E-K-1213}
As $N \to \infty$ and $h\to 0$, it holds that
\begin{align*}
E[K_{12}(\alphabh[])K_{13}(\alphabh[])]
& \sim h^2 E[f_{\eps[]}^2(y - \Xb[1]\alphab[])].
\end{align*}
\begin{proof}
Setting $\Ib = \{1,2,3\}$, we have that at leading-order, 
\begin{align*}
E[K_{12}(\alphabh[])K_{13}(\alphabh[])]
\sim E[K_{12}(\alphabh[I])K_{13}(\alphabh[I])]
& = \int \cdots \int K_h(\yt(\alphabh[I],\xb[1],\xb[2]) - e_2)K_h(\yt(\alphabh[I],\xb[1],\xb[3]) - e_3) \\
& \times f_{\eps[2]}(e_2)f_{\eps[3]}(e_3) f_{\epsb[23-]}(\eb[23-]) f_{\Xb[]}(\xb[]) \ d\xb[] d \eb[23-] de_2 de_3. 
\end{align*}
Performing the change of variables $-r_1 = \yt(\alphabh[I],\xb[1],\xb[2]) - e_2)/h$ and $-r_2 = \yt(\alphabh[I],\xb[1],\xb[3]) - e_3)/h$, we get
\begin{align*}
E[K_{12}(\alphabh[])K_{13}(\alphabh[])]
& \sim h^2 \int \cdots \int K(r_1)K(r_2) f_{\eps[2]}(\yt(\alphabh[I],\xb[1],\xb[2]) + hr_1)f_{\eps[3]}(\yt(\alphabh[I],\xb[1],\xb[3]) + hr_2) \\
& \times f_{\epsb[23-]}(\eb[23-]) f_{\Xb[]}(\xb[]) \ d\xb[] d \eb[23-] dr_1 dr_2. 
\end{align*}
Taylor expanding with respect to $h$ about $0$, and using \eqref{eq:int-K}, gives
\begin{align*}
E[K_{12}(\alphabh[])K_{13}(\alphabh[])]
& \sim h^2 \int \cdots \int f_{\eps[2]}(\yt(\alphabh[I],\xb[1],\xb[2]))f_{\eps[3]}(\yt(\alphabh[I],\xb[1],\xb[3])) f_{\epsb[23-]}(\eb[23-]) f_{\Xb[]}(\xb[]) \ d\xb[] d \eb[23-] \\
& = h^2 E[f_{\eps[]}(\yt(\alphabh[I],\Xb[1],\Xb[2]))f_{\eps[]}(\yt(\alphabh[I],\Xb[1],\Xb[3]))].
\end{align*}
Finally, taking the leading-order approximation $\alphab[]$ of $\alphabh[I]$, we have that
\begin{align*}
E[K_{12}(\alphabh[])K_{13}(\alphabh[])]
& \sim h^2 E[f_{\eps[]}(y - \Xb[1]\alphab[])f_{\eps[]}(y - \Xb[1]\alphab[])]
= h^2 E[f_{\eps[]}^2(y - \Xb[1]\alphab[])].
\end{align*}
\end{proof}
\end{lemma}

%
%
%
%
%
%
%
%
\begin{lemma} \label{lem:E-K-1232}
As $N \to \infty$ and $h\to 0$, it holds that
\begin{align*}
E[K_{12}(\alphabh[])K_{32}(\alphabh[])]
\sim h^2 \int_{R} f_{\eps[]}(y - \xb[1]^T \alphab[]) f_{\Xb[1]}(\xb[1]) f_{\Xb[3]}(\xb[3]) \ d\xb[1] d\xb[3].
\end{align*}
where the region of integration is $R =\{(\xb[1],\xb[3]) : (\xb[1] - \xb[3])^T \alphab[]=0\}$.
\begin{proof}
We have $E[K_{12}(\alphabh[])K_{32}(\alphabh[])] = E[K_{12}(\alphab[])K_{32}(\alphab[])] + (E[K_{12}(\alphabh[])K_{32}(\alphabh[])]-E[K_{12}(\alphab[])K_{32}(\alphab[])])$, which corresponds to the decomposition approach used in \cite{stove2012convolution}. In that work, it was established that the first term dominates asymptotically so it suffices to consider
\begin{align*}
E[K_{12}(\alphab[])K_{32}(\alphab[])]
& \sim \int \cdots \int K_h(y - \xb[1]^T \alphab[] - e_2)K_h(y - \xb[3]^T \alphab[] - e_2) f_{\eps[2]}(e_2) f_{\Xb[1]}(\xb[1]) f_{\Xb[3]}(\xb[3]) \ d\xb[1] d\xb[3] de_2. 
\end{align*}
Performing the change of variables $-r_1 = (y - \xb[1]^T \alphab[] - e_2)/h$, using the fact that $K$ is a symmetric function, and taking the leading-order Taylor approximation of $f_{\eps[2]}(y - \xb[1]^T \alphab[] + hr_1)$, we find that
\begin{align*}
E[K_{12}(\alphab[])K_{32}(\alphab[])]
& = h \int \cdots \int K(r_1)K\bigg(\frac{(\xb[1] - \xb[3])^T \alphab[]}{h} - r_1\bigg) f_{\eps[2]}(y - \xb[1]^T \alphab[] + hr_1) f_{\Xb[1]}(\xb[1]) f_{\Xb[3]}(\xb[3]) \ d\xb[1] d\xb[3] dr_1 \\
& \sim h \int \cdots \int K(r_1)K\bigg(\frac{(\xb[1] - \xb[3])^T \alphab[]}{h} - r_1\bigg) f_{\eps[2]}(y - \xb[1]^T \alphab[]) f_{\Xb[1]}(\xb[1]) f_{\Xb[3]}(\xb[3]) \ d\xb[1] d\xb[3] dr_1 \\
& = h^2 \int \cdots \int h^{-1} K^*\bigg(\frac{(\xb[1] - \xb[3])^T \alphab[]}{h}\bigg) f_{\eps[2]}(y - \xb[1]^T \alphab[]) f_{\Xb[1]}(\xb[1]) f_{\Xb[3]}(\xb[3]) \ d\xb[1] d\xb[3] \\
& \sim h^2 \int_{R} f_{\eps[]}(y - \xb[1]^T \alphab[]) f_{\Xb[1]}(\xb[1]) f_{\Xb[3]}(\xb[3]) \ d\xb[1] d\xb[3],
\end{align*}
as $h \to 0$, up to a constant factor that is irrelevant to the asymptotic analysis, where $K^*(a)$ is the Gaussian function given by the convolution $K^*(a) = \int K(b)K(a-b)db$. See \cite[Supp. Material, Proof of Theorem 3] {stove2012convolution} for the analogous approach in the case of the Nadaraya-Watson-enhanced convolution estimator; in particular, the expression above has a correspondence to the second term in \cite[(19)] {stove2012convolution}.
\end{proof}
\end{lemma}

\begin{lemma} \label{lem:E-feps-2-feps-4}
As $N \to \infty$, it holds that
\begin{align*}
E[f_{\eps[]}(\yt(\alphabh[I],\Xb[1],\Xb[2]))f_{\eps[]}(\yt(\alphabh[I],\Xb[3],\Xb[4]))]
& \sim \fY^2(y)
+ N^{-1} \sigmaeps^2 \sum_{p_1,p_2=0}^J E[(\Phii[N]^{-1} \Xb[1])_{p_1}(\Phii [N]^{-1} \Xb[1])_{p_2}] \\
& \times (\fY(y) E[f_{\eps[]}''(y - \Xb[1]^T \alphab[]) (\Xb[2] - \Xb[1])_{p_1}(\Xb[2] - \Xb[1])_{p_2}]  \\
& + E[f_{\eps[]}'(y - \Xb[1]^T \alphab[]) (\Xb[2] - \Xb[1])_{p_1}]E[f_{\eps[]}'(y - \Xb[3]^T \alphab[]) (\Xb[4] - \Xb[3])_{p_2}]).
\end{align*}
\begin{proof}
Taylor expanding with respect to $\alphabh[I]$ at about $\alphab[]$, and taking the expectation, we get
\begin{align*}
& E[f_{\eps[]}(\yt(\alphabh[I],\Xb[1],\Xb[2]))f_{\eps[]}(\yt(\alphabh[I],\Xb[3],\Xb[4]))] \\
&\sim E[f_{\eps[]}(y - \Xb[1]^T \alphab[])] E[f_{\eps[]}(y - \Xb[3]^T \alphab[])] \\
&+ \sum_{p_1=0}^J E[(\alphabh[I]-\alphab[])_{p_1}] E\bigg[\frac{\partial f_{\eps[]}(\yt(\alphabh[],\Xb[1],\Xb[2]))f_{\eps[]}(\yt(\alphabh[I],\Xb[3],\Xb[4]))}{\partial(\alphabh[I])_{p_1}} \bigg|_{\alphabh[I]=\alphab[]}\bigg] \\
& + \frac{1}{2} \sum_{p_1,p_2=0}^J E[(\alphabh[I]-\alphab[])_{p_1}(\alphabh[I]-\alphab[])_{p_2}] E\bigg[\frac{\partial^2 f_{\eps[]}(\yt(\alphabh[I],\Xb[1],\Xb[2]))f_{\eps[]}(\yt(\alphabh[I],\Xb[3],\Xb[4]))}{\partial(\alphabh[I])_{p_1}\partial(\alphabh[I])_{p_2}} \bigg|_{\alphabh[I]=\alphab[]}\bigg],
%
\end{align*}
where we have ignored the remainder since it can be neglected as shown in \eqref{lem:E-feps-2}. Now, $E[f_{\eps[]}(y - \Xb[1]^T \alphab[])] E[f_{\eps[]}(y - \Xb[3]^T \alphab[])] = \fY^2(y)$. The first order term vanishes since $E[(\alphabh[I] - \alphab[])p_1] = 0$. It was shown in \eqref{lem:E-feps-2} that $E[(\alphabh[I]-\alphab[])_{p_1}(\alphabh[I]-\alphab[])_{p_2}] = N^{-1} \sigmaeps^2 E[(\Phii[N]^{-1} \Xb[1])_{p_1}(\Phii [N]^{-1} \Xb[1])_{p_2}]$. Next,
\begin{align*}
\frac{\partial^2 f_{\eps[]}(\yt(\alphabh[I],\Xb[1],\Xb[2]))f_{\eps[]}(\yt(\alphabh[I],\Xb[3],\Xb[4]))}{\partial(\alphabh[I])_{p_1}\partial(\alphabh[I])_{p_2}}\bigg|_{\alphabh[I]=\alphab[]}
&= f_{\eps[]}''(y - \Xb[1]^T \alphab[]) f_{\eps[]}(y - \Xb[3]^T \alphab[]) (\Xb[2] - \Xb[1])_{p_1}(\Xb[2] - \Xb[1])_{p_2}  \\
& + f_{\eps[]}'(y - \Xb[1]^T \alphab[]) f_{\eps[]}'(y - \Xb[3]^T \alphab[]) (\Xb[2] - \Xb[1])_{p_1}(\Xb[4] - \Xb[3])_{p_2} \\
& + f_{\eps[]}'(y - \Xb[1]^T \alphab[]) f_{\eps[]}'(y - \Xb[3]^T \alphab[]) (\Xb[2] - \Xb[1])_{p_2}(\Xb[4] - \Xb[3])_{p_1} \\
& + f_{\eps[]}(y - \Xb[1]^T \alphab[]) f_{\eps[]}''(y - \Xb[3]^T \alphab[]) (\Xb[4] - \Xb[3])_{p_1}(\Xb[4] - \Xb[3])_{p_2}.
\end{align*}
Therefore
\begin{align*}
& E\bigg[\frac{\partial^2 f_{\eps[]}(\yt(\alphabh[I],\Xb[1],\Xb[2]))f_{\eps[]}(\yt(\alphabh[I],\Xb[3],\Xb[4]))}{\partial(\alphabh[I])_{p_1}\partial(\alphabh[I])_{p_2}}\bigg|_{\alphabh[I]=\alphab[]}\bigg] \\
&= 2 (E[f_{\eps[]}''(y - \Xb[1]^T \alphab[]) f_{\eps[]}(y - \Xb[3]^T \alphab[]) (\Xb[2] - \Xb[1])_{p_1}(\Xb[2] - \Xb[1])_{p_2}]  \\
& + E[f_{\eps[]}'(y - \Xb[1]^T \alphab[]) f_{\eps[]}'(y - \Xb[3]^T \alphab[]) (\Xb[2] - \Xb[1])_{p_1}(\Xb[4] - \Xb[3])_{p_2}]) \\
&= 2 (\fY(y) E[f_{\eps[]}''(y - \Xb[1]^T \alphab[]) (\Xb[2] - \Xb[1])_{p_1}(\Xb[2] - \Xb[1])_{p_2}]  \\
& + E[f_{\eps[]}'(y - \Xb[1]^T \alphab[]) (\Xb[2] - \Xb[1])_{p_1}]E[f_{\eps[]}'(y - \Xb[3]^T \alphab[]) (\Xb[4] - \Xb[3])_{p_2}]).
\end{align*}
\end{proof}
\end{lemma}

%
%
%
%
%
%
%
%
\begin{lemma} \label{lem:E-K-1234}
As $N \to \infty$ and $h\to 0$, it holds that
\begin{align*}
E[K_{12}(\alphabh[])K_{34}(\alphabh[])]
& \sim h^2 \fY^2(y)
+ h^4 \mu_K^2 \fY(y) \fY''(y) \\
& + h^2 N^{-1} \sigmaeps^2 \sum_{p_1,p_2=0}^J E[(\Phii[N]^{-1} \Xb[1])_{p_1}(\Phii[N]^{-1} \Xb[1])_{p_2}] \\
& \times (\fY(y) E[f_{\eps[]}''(y - \Xb[1]^T \alphab[]) (\Xb[2] - \Xb[1])_{p_1}(\Xb[2] - \Xb[1])_{p_2}] \\
& + E[f_{\eps[]}'(y - \Xb[1]^T \alphab[]) (\Xb[2] - \Xb[1])_{p_1}]E[f_{\eps[]}'(y - \Xb[3]^T \alphab[]) (\Xb[4] - \Xb[3])_{p_2}]).
\end{align*}
\begin{proof}
Setting $\Ib = \{1,2,3,4\}$, we have that at leading-order,
\begin{align*}
E[K_{12}(\alphabh[])K_{34}(\alphabh[])] \sim E[K_{12}(\alphabh[I])K_{34}(\alphabh[I])]
& = \int \cdots \int K_h(\yt(\alphabh[I],\xb[1],\xb[2]) - e_2)K_h(\yt(\alphabh[I],\xb[3],\xb[4]) - e_4) \\
& \times f_{\eps[2]}(e_2) f_{\eps[4]}(e_4) f_{\epsb[1234-]}(\eb[1234-]) f_{\Xb[]}(\xb[]) \ d\xb[] d\eb[1234-] de_2 de_4. 
\end{align*}
Performing the change of variables $-r_1 = \yt(\alphabh[I],\xb[1],\xb[2]) - e_2)/h$ and $-r_2 = \yt(\alphabh[I],\xb[3],\xb[4]) - e_4)/h$, we get
\begin{align*}
E[K_{12}(\alphabh[])K_{34}(\alphabh[])]
& \sim h^2 \int \cdots \int K(r_1)K(r_2) f_{\eps[2]}(\yt(\alphabh[I],\xb[1],\xb[2]) + hr_1) f_{\eps[4]}(\yt(\alphabh[I],\xb[3],\xb[4]) + hr_2) \\
& \times f_{\epsb[1234-]}(\eb[1234-]) f_{\Xb[]}(\xb[]) \ d\xb[] d\eb[1234-] dr_1 dr_2. 
\end{align*}
Taylor expanding with respect to $h$ about $0$, and using \eqref{eq:int-K} and \eqref{eq:int-K-defs}, we get
\begin{align*}
E[K_{12}(\alphabh[])K_{34}(\alphabh[])]
& \sim h^2 \int \cdots \int f_{\eps[2]}(\yt(\alphabh[I],\xb[1],\xb[2]))f_{\eps[4]}(\yt(\alphabh[I],\xb[3],\xb[4])) f_{\epsb[1234-]}(\eb[1234-]) f_{\Xb[]}(\xb[]) \ d\xb[] d\eb[1234-] \\
&+ h^4 \frac{\mu_K^2}{2} \int \cdots \int (f_{\eps[2]}(\yt(\alphabh[I],\xb[1],\xb[2]))f_{\eps[4]}''(\yt(\alphabh[I],\xb[3],\xb[4])) 
+ f_{\eps[2]}''(\yt(\alphabh[I],\xb[1],\xb[2])) f_{\eps[4]}(\yt(\alphabh[I],\xb[3],\xb[4]))) \\
& \times f_{\epsb[1234-]}(\eb[1234-]) f_{\Xb[]}(\xb[]) \ d\xb[] d\eb[1234-] \\ 
%
& = h^2 E[f_{\eps[]}(\yt(\alphabh[I],\Xb[1],\Xb[2]))f_{\eps[]}(\yt(\alphabh[I],\Xb[3],\Xb[4]))] \\
& + h^4 \frac{\mu_K^2}{2}(E[f_{\eps[]}(\yt(\alphabh[I],\Xb[1],\Xb[2]))f_{\eps[]}''(\yt(\alphabh[I],\Xb[3],\Xb[4]))] + E[f_{\eps[]}''(\yt(\alphabh[I],\Xb[1],\Xb[2]))f_{\eps[]}''(\yt(\alphabh[I],\Xb[3],\Xb[4]))]) \\
& = h^2 E[f_{\eps[]}(\yt(\alphabh[I],\Xb[1],\Xb[2]))f_{\eps[]}(\yt(\alphabh[I],\Xb[3],\Xb[4]))]
+ h^4 \mu_K^2 E[f_{\eps[]}(\yt(\alphabh[I],\Xb[1],\Xb[2]))f_{\eps[]}''(\yt(\alphabh[I],\Xb[3],\Xb[4]))].
\end{align*}
By Lemma \ref{lem:E-feps-2-feps-4}, we have that at leading-order, 
\begin{align*}
E[f_{\eps[]}(\yt(\alphabh[I],\Xb[1],\Xb[2]))f_{\eps[]}(\yt(\alphabh[I],\Xb[3],\Xb[4]))]
& \sim \fY^2(y)
+ N^{-1} \sigmaeps^2 \sum_{p_1,p_2=0}^J E[(\Phii[N]^{-1} \Xb[1])_{p_1}(\Phii[N]^{-1} \Xb[1])_{p_2}] \\
& \times (\fY(y) E[f_{\eps[]}''(y - \Xb[1]^T \alphab[]) (\Xb[2] - \Xb[1])_{p_1}(\Xb[2] - \Xb[1])_{p_2}] \\
& + E[f_{\eps[]}'(y - \Xb[1]^T \alphab[]) (\Xb[2] - \Xb[1])_{p_1}]E[f_{\eps[]}'(y - \Xb[3]^T \alphab[]) (\Xb[4] - \Xb[3])_{p_2}]).
%
\end{align*}
Finally, taking the leading-order approximation $\alphab[]$ of $\alphabh[I]$, and using \eqref{eq:fYk-conv-form}, we have
\begin{align*}
E[f_{\eps[]}(\yt(\alphabh[I],\Xb[1],\Xb[2]))f_{\eps[]}''(\yt(\alphabh[I],\Xb[3],\Xb[4]))]
\sim E[f_{\eps[]}(y - \Xb[1] \alphab[])]E[f_{\eps[]}''(y - \Xb[3] \alphab[])]
= \fY(y) \fY''(y).
\end{align*}

\end{proof}
\end{lemma}
\section{} \label{sec:appendix-E-PhiiN-inv-sqr-ord-mag}
For convenience, throughout this section we often use the notation $\lesssim$ to represent an expression that holds up an asymptotically irrelevant constant factor as $N \to \infty$. For example, instead of writing $5N^{-1} + 3N^{-2} \le C N^{-1}$ for $C>0$ as $N \to \infty$, we write $5N^{-1} + 3N^{-2} \lesssim N^{-1}$ as $N \to \infty$. Before proving Lemma \ref{lem:E-prod-PhiiN-Xi}, which is the primary objective of this section, we need several lemmas.

%
%
%
%
%
%
%
%
\begin{lemma} \label{lem:E-A-v}
Let $\Ab[] \in \mathbb{R}^{(J+1)\times (J+1)}$ be a matrix such that $E[|(\Ab[])_{ij}|^r]$ and $E[|(\Ab[])_{kl}|^r]$ have the same scaling with respect to $N$, and let $\vb[] \in \mathbb{R}^{J+1}$ be a vector such that $E[|(\vb[])_{j}|^r]$ and $E[|(\vb[])_{l}|^r]$ have the same scaling with respect to $N$, as $N \to \infty$, for $r \ge 1$ and $i,j,k,l \in \{0,1,\dots,J\}$. In particular, assume that
\begin{align}
E[(\Ab[])_{ij}^4] & \lesssim E[(\Ab[])_{11}^4], \label{eq:A-ij-le} \\
E[(\vb[])_{j}^4] & \lesssim E[(\vb[])_{1}^4], \label{eq:v-j-le}
\end{align}
as $N \to \infty$, for $i,j \in \{0,1,\dots,J\}$. Then, as $N \to \infty$,
\begin{align*}
E[(\Ab[]\vb[])_i^2] \lesssim E[(\Ab[])_{11}^4]^{1/2} E[(\vb[])_{1}^4]^{1/2}.
\end{align*}
\begin{proof}
The $i$-th element of the vector $\Ab[]\vb[]$ is $(\Ab[]\vb[])_i = \sum_{j=0}^J (\Ab[])_{ij} (\vb[])_{j}$. Therefore,
\begin{align*}
E[(\Ab[]\vb[])_i^2]
= E\bigg[\bigg(\sum_{j=0}^J (\Ab[])_{ij} (\vb[])_{j}\bigg)^2\bigg]
= \sum_{j=0}^J \sum_{k=0}^J E[(\Ab[])_{ij} (\Ab[])_{ik} (\vb[])_{j} (\vb[])_{k}].
\end{align*}
Then, we use the Cauchy-Schwarz inequality, followed by \eqref{eq:A-ij-le} and \eqref{eq:v-j-le}, to get
\begin{align*}
E[(\Ab[]\vb[])_i^2]
& \le \sum_{j=0}^J \sum_{k=0}^J E[((\Ab[])_{ij} (\Ab[])_{ik})^2]^{1/2} E[((\vb[])_{j} (\vb[])_{k})^2]^{1/2} \\
& \le \sum_{j=0}^J \sum_{k=0}^J E[(\Ab[])_{ij}^4]^{1/4} E[(\Ab[])_{ik}^4]^{1/4} E[(\vb[])_{j}^4]^{1/4} E[(\vb[])_{k}^4]^{1/4} \\
& \lesssim (J+1)^2 E[(\Ab[])_{11}^4]^{1/2} E[(\vb[])_{1}^4]^{1/2} \\
& \lesssim E[(\Ab[])_{11}^4]^{1/2} E[(\vb[])_{1}^4]^{1/2}.
\end{align*}
\end{proof}
\end{lemma}


%
%
%
%
%
%
%
%
\begin{proposition} \label{prop:E-S_N}
Let $Z_1,Z_2,\dots,$ be i.i.d positive random variables, and let $S_N = \sum_{i=1}^N Z_i$. Suppose $E[Z_1^r] < \infty$ for $r \ge 2$,. Then, as $N \to \infty$, $E[|S_N|^r]$ scales as
\begin{align*}
E[|S_N|^r] \sim C N^r,
\end{align*}
for $C > 0$.
\begin{proof}
Since $Z_i$ is positive, there exists $C_1 > 0$ such that $Z_i \ge C_1$. Then $(\sum_{i=1}^N Z_i)^r \ge C_1^r N^r=C_2N^r$, and so $E[|S_N|^r] \ge C_2N^r$.

Denote by $\mu_Z = E[Z_i]$ and $\sigma_Z^2 = \Var[Z_i]$, and note that by Minkowski's inequality,
\begin{align*}
E\bigg[\bigg|\frac{S_N}{\sigma_Z \sqrt{N}}\bigg|^r\bigg]^{1/r}
& = E\bigg[\bigg|\frac{S_N-N\mu_Z}{\sigma_Z \sqrt{N}} + \frac{N\mu_Z}{\sigma_Z \sqrt{N}}\bigg|^r\bigg]^{1/r} \\
& \le E\bigg[\bigg|\frac{S_N-N\mu_Z}{\sigma_Z \sqrt{N}}\bigg|^r\bigg]^{1/r} + \frac{\sqrt{N}\mu_Z}{\sigma_Z}.
\end{align*}
Then, by Theorem \cite[Thm. 5.1]{gut2013probability}, there exists $C_3 > 0$ such that as $N \to \infty$,
\begin{align*}
E\bigg[\bigg|\frac{S_N-N\mu_Z}{\sigma_Z \sqrt{N}}\bigg|^r\bigg]^{1/r} \le C_3 E[|W|^r]^{1/r},
\end{align*}
where $W \sim \mathcal{N}(0,1)$. This means that means that as $N \to \infty$,
\begin{align*}
E\bigg[\bigg|\frac{S_N}{\sigma_Z \sqrt{N}}\bigg|^r\bigg]^{1/r}
\le C_3 E[|W|^r]^{1/r} + \frac{\sqrt{N}\mu_Z}{\sigma_Z}
\le C_4 \sqrt{N}.
\end{align*}
Therefore,
\begin{align*}
E[|S_N|^r]
& \le C_5 N^r.
\end{align*}
Since $C_2 N^r \le E[S_N] \le C_5 N^r$, we have $E[|S_N|^r] \sim C N^r$.
\end{proof}
\end{proposition}


%
%
%
%
%
%
%
%
\begin{lemma} \label{lem:E-Det-r}
Let $r \ge 1$. As $N \to \infty$,
\begin{align*}
E[\Det^{r}(\Xb[]^T\Xb[])] \sim C N^{r(J+1)},
\end{align*}
for $C > 0$.
\begin{proof}
The formula for the determinant of a matrix can be written as a summation that includes a term given by the product of the diagonals of the matrix. It suffices to consider only this term since the scaling with respect to $N$ is common to all of the terms in summation. Therefore, since the diagonal of $\Xb[]^T \Xb[]$ is $[N,\sum_{i=1}^N \X[i1]^2,\sum_{i=1}^N \X[i2]^2\dots,\sum_{i=1}^N \X[iJ]^2]^T
$, we have that as $N \to \infty$, the determinant scales as
\begin{align*}
E[\Det^r(\Xb[]^T \Xb[])]
\sim C_1 E\bigg[\bigg(N \prod_{j=1}^J \sum_{i=1}^N \X[ij]^2\bigg)^r\bigg]
= C_1 N^r \prod_{j=1}^J E\bigg[\bigg(\sum_{i=1}^N \X[ij]^2\bigg)^r\bigg],
\end{align*}
where we used the fact that $\sum_{i=1}^N \X[ij]^2$ is independent of $\sum_{i=1}^N \X[ik]^2$ for $j \neq k$ in the last line. Then, by Proposition \ref{prop:E-S_N}, $E[(\sum_{i=1}^N \X[ij]^2)^r] \sim C_2 N^r$, and thus
\begin{align*}
E[\Det^r(\Xb[]^T \Xb[])]
\sim C_3 N^r \prod_{j=1}^J N^r
= C_3 N^{r(J+1)}.
\end{align*}
\end{proof}
\end{lemma}

%
%
%
%
%
%
%
%
\begin{corollary} \label{cor:E-Det-neg-r}
Let $r \ge 1$. As $N \to \infty$,
\begin{align*}
E[\Det^{-r}(\Xb[]^T \Xb[])]
\sim C N^{-r(J+1)}.
\end{align*}
\begin{proof}
Define $g(x) = x^{-1}$, and note that $\Det^r(\Xb[]^T \Xb[])$, is a non-negative random variable. Taylor expanding $g(\Det^r(\Xb[]^T \Xb[]))$ at $\Det^r(\Xb[]^T \Xb[]) = E[\Det^r(\Xb[]^T \Xb[])]$ and then taking the leading-order approximation of its expectation, we have that
\begin{align*}
E[\Det^{-r}(\Xb[]^T \Xb[])]
= E[g(\Det^r(\Xb[]^T \Xb[]))]
\sim g(E[\Det^{r}(\Xb[]^T \Xb[])])
= E[\Det^{r}(\Xb[]^T \Xb[])]^{-1},
\end{align*}
as $N \to \infty$.
The result follows by Lemma \ref{lem:E-Det-r}.
\end{proof}
\end{corollary}

%
%
%
%
%
%
%
%
\begin{corollary} \label{cor:E-Adj-r}
Let $r,s \ge 1$. As $N \to \infty$,
\begin{align*}
E[(\Adj^r(\Xb[]^T\Xb[])_{ij})^s] \sim C N^{rsJ},
\end{align*}
for $C > 0$ and $i,j \in \{0,1,\dots,J\}$.
\begin{proof}
The adjugate matrix of $\Xb[]^T\Xb[]$ is the transpose of its cofactor matrix. The elements of the cofactor matrix are themselves determinants of $J\times J$ submatrices formed from the elements of the $(J+1)\times(J+1)$-dimensional matrix $\Xb[]^T\Xb[]$. In particular, the scaling derived in Lemma \ref{lem:E-Det-r} applies to these determinants. Denote by $\Xb[J-]$ the matrix $\Xb[]$ in the multiple regression model \eqref{eq:samnple-mult-reg} with the $J$-th column removed. Then $\Adj^r(\Xb[]^T\Xb[])_{ij}$ scales like $\Det^r(\Xb[J-]^T\Xb[J-])$ with respect to $N$, which in turn means that $(\Adj^r(\Xb[]^T\Xb[])_{ij})^s$ scales like $\Det^{rs}(\Xb[J-]^T\Xb[J-])$. To be precise,
\begin{align*}
E[\Adj^r(\Xb[]^T\Xb[])_{ij}]
\sim C_1 E[\Det^{sr}(\Xb[J-]^T\Xb[J-])]
\sim C_2 N^{rs J},
\end{align*}
where we used Lemma \eqref{lem:E-Det-r} with the $(J+1)\times(J+1)$-dimensional matrix $\Xb[]^T\Xb[]$ replaced by the $J\times J$-dimensional matrix $\Xb[J-]^T\Xb[J-]$ to obtain the final expression.
\end{proof}
\end{corollary}

%
%
%
%
%
%
%
%
\begin{corollary} \label{cor:E-Det-Adj-r}
Let $r,s \ge 1$. As $N \to \infty$,
\begin{align*}
E[((\Xb[]^T \Xb[])^{-r})_{11}^s] \lesssim N^{-rs}.
\end{align*}
\begin{proof}
Since $(\Xb[]^T \Xb[])^{-r} = \Det^{-r}(\Xb[]^T \Xb[])\Adj^r(\Xb[]^T \Xb[])$, using properties of the determinant and adjugate matrix, followed by the Cauchy-Schwarz inequality, we get
\begin{align*}
E[((\Xb[]^T \Xb[])^{-r})_{11}^s]
& = E[(\Det^{-r}(\Xb[]^T \Xb[]) \Adj^r(\Xb[]^T \Xb[]))_{11}^s] \\
& = E[\Det^{-sr}(\Xb[]^T \Xb[]) (\Adj^r(\Xb[]^T \Xb[]))_{11}^s] \\
& \le E[(\Det^{-2sr}(\Xb[]^T \Xb[])]^{1/2} E[\Adj^r(\Xb[]^T \Xb[]))_{11}^{2s}]^{1/2}.
\end{align*}
Now, by Corollary \eqref{cor:E-Det-neg-r}, $E[(\Det^{-2sr}(\Xb[]^T \Xb[])]^{1/2} \sim N^{sr(J+1)}$, and by 
Corollary \eqref{cor:E-Adj-r}, $E[\Adj^r(\Xb[]^T \Xb[]))_{11}^{2s}]^{1/2} \sim N^{rsJ}$. Hence,
\begin{align*}
E[((\Xb[]^T \Xb[])^{-r})_{11}^s]
\lesssim N^{-rs(J+1)} N^{rsJ}
=  N^{-rs}.
\end{align*}
\end{proof}
\end{corollary}

%
%
%
%
%
%
%
%
\begin{lemma} \label{lem:X-T-X-inv-X-1-ineq}
As $N \to \infty$, it holds that
\begin{align*}
E[((\Xb[]^T \Xb[])^{-1}\Xb[1])_{i}^k] = O(N^{-k}).
\end{align*}
\begin{proof}
By Lemma \ref{lem:E-A-v},
\begin{align*}
E[((\Xb[]^T \Xb[])^{-1}\Xb[1])_{i}^k]
\lesssim E[((\Xb[]^T \Xb[])^{-1})_{11}^{2k}]^{1/2} E[(\Xb[1])_{1}^{2k}]^{1/2}
\lesssim E[((\Xb[]^T \Xb[])^{-1})_{11}^{2k}]^{1/2},
\end{align*}
since $\Xb[1]$ is independent of $N$. Then applying \eqref{cor:E-Det-Adj-r} with $(r,s) = (1,2k)$, we get
\begin{align*}
E[((\Xb[]^T \Xb[])^{-1})_{11}^{2k}]^{1/2}
\lesssim 
(N^{-2k})^{1/2}
= N^{-k}.
\end{align*}
\end{proof}
\end{lemma}

%
%
%
%
%
%
%
%
\begin{lemma} \label{lem:E-prod-PhiiN-Xi}
As $N\to \infty$, it holds that
\begin{align*}
E\bigg[\prod_{i=1}^k (\Phii[N]^{-1} \Xb[1])_{p_i}\bigg] = O(1),
\end{align*}
for $k \in \mathbb{N}$.
\begin{proof}
Since all the elements of the vector $\Phii[N]^{-1} \Xb[1]$ have the same order of magnitude with respect to $N$,
\begin{align*}
E\bigg[\prod_{i=1}^k (\Phii[N]^{-1} \Xb[1])_{p_i}\bigg]
\lesssim E[(\Phii[N]^{-1} \Xb[1])_{p_1}^k].
\end{align*}
Therefore, recalling \eqref{eq:Phi}, and then using \eqref{lem:X-T-X-inv-X-1-ineq}, we have that
\begin{align*}
E\bigg[\prod_{i=1}^k (\Phii[N]^{-1} \Xb[1])_{p_i}\bigg]
\lesssim E[(\Phii[N]^{-1} \Xb[1])_{p_1}^k]
= N^k E[((\Xb[]^T \Xb[])^{-1} \Xb[1])_{p_1}^k]
\lesssim N^k N^{-k}
= 1.
\end{align*}
\end{proof}
\end{lemma}

\end{appendices}

\end{document}